\documentstyle[12pt]{article}
\font\teneufm=eufm10
\font\seveneufm=eufm7
\font\fiveeufm=eufm5
\newfam\frakturfam

\textfont\frakturfam=\teneufm
\scriptfont\frakturfam=\seveneufm
\scriptscriptfont\frakturfam=\fiveeufm

\newtheorem{pr}{Proposition}[section]
\newtheorem{df}{Definition}[section]
\newtheorem{lm}{Lemma}[section]
\newtheorem{theor}{Theorem}[section]
\newtheorem{co}{Corollary}[section]
\newtheorem{rem}{Remark}[section]
\newtheorem{prob}{Problem}[section]
\def\bee{\begin{eqnarray}}
\def\bes{\begin{eqnarray*}}
\def\eee{\end{eqnarray}}
\def\ees{\end{eqnarray*}}
\def\le{{\langle}}
\def\re{{\rangle}}

\def\a{\alpha}
\def\b{\beta}
\def\te{\theta}
\def\g{\gamma}

\def\s{\sigma}
\def\t{\tau}
\def\d{\partial}
\def\l{\lambda}

\def\Proof{{\sl Proof.}\ }

\title{Defining relations of the tame automorphism group of polynomial algebras in three variables}

\begin{document}
\date{}
\maketitle

\begin{center}

{\bf U.\,U.\,Umirbaev}\\
Eurasian National University,\\
 Astana, 010008, Kazakhstan \\
e-mail: {\em umirbaev@yahoo.com}

\end{center}

\begin{abstract}  
We describe a set of defining relations of the tame automorphism group $TA_3(F)$  
of the polynomial algebra $F[x_1,x_2,x_3]$ in variables $x_1,x_2,x_3$ over 
an arbitrary field $F$ of characteristic $0$.  
\end{abstract}



\section{Introduction}

\hspace*{\parindent}

Let $A_n=F[x_1,x_2,\ldots,x_n]$ be the polynomial algebra 
in the variables $x_1,x_2,\ldots,x_n$ over a field $F$, and 
let  $GA_n(F)=Aut\,A_n$ be the automorphism group of $A_n$.
  Let  $\phi = (f_1,f_2,\ldots,f_n)$ denote an automorphism 
  $\phi$ of $A_n$ such that 
   $\phi(x_i)=f_i,\, 1\leq i\leq n$. An automorphism  
\bes
\s(i,\a,f) = (x_1,\ldots,x_{i-1}, \a x_i+f, x_{i+1},\ldots,x_n),
\ees
where $0\neq\a\in F,\ f\in F[x_1,\ldots,x_{i-1},x_{i+1},\ldots,x_n]$, 
 is called {\em elementary}. The subgroup $TA_n(F)$ of $GA_n(F)$ generated by all 
  elementary automorphisms is called the {\em tame automorphism group}, 
 and the elements of this subgroup are called the {\em tame automorphisms} 
 of $A_n$. Nontame automorphisms of $A_n$ are called {\em wild}.

\smallskip

 It is well known \cite{Czer,Jung,Kulk,Makar} that the automorphisms 
of polynomial algebras and free associative algebras in two variables are tame. 
It was recently proved in    
\cite{SU4,Umi24,Umi25,SU3} that the well-known Nagata automorphism (see \cite{Nagata})
\bes
\sigma&=&(x+(x^2-yz)z,\,y+2(x^2-yz)x+(x^2-yz)^2z,\,z)
\ees
of the polynomial algebra $F[x,y,z]$ over a field $F$ of characteristic $0$ is wild.

\smallskip

It is well known (see, for example \cite{Cohn}) that the groups of automorphisms 
of polynomial algebras and free associative algebras in two variables are isomorphic  
and have a nice representation as a free product of groups. Similar results are well 
known for two dimensional Cremona groups \cite{Gizat, Iskovskikh, Wright1}. 
In combinatorial group theory (see, for example \cite{Vogt}) there are several well-known descriptions 
of the group of automorphisms of free groups by generators and defining relations. 
Main part of the  
 investigations of the automorphism groups in the case of free linear 
algebras consist in finding a system of generators. In 1968 P.Cohn \cite{Cohn2} proved that the automorphisms of free 
Lie algebras with a finite set of generators are tame. 

\smallskip

In this paper we describe a set of defining relations of the tame automorphism group $TA_3(F)$ of the polynomial 
algebra $A=F[x_1,x_2,x_3]$ over a field $F$ of characteristic $0$.   
The discovered set of defining relations (\ref{f3.2})--(\ref{f3.4}) is new and very natural for free linear algebras. 
This opens the way to the representation theory of the automorphism groups of free associative algebras, free Lie algebras and the others. Using the obtained presentation of $TA_3(F)$ we proved (see \cite{UUU}) that the well-known Anick automorphism    of the free associative algebra  in three variables is wild. 

The paper is organized as follows. 
In Section~2 we give some preliminaries from \cite{Umi24,Umi25}. 
It was proved  in \cite{Umi25} that every tame automorphism of the polynomial algebra  
$A$ of degree greater than $3$ admits either an elementary reduction or a reduction of types I--IV. 
In this paper we use some extension of the definition of automorphisms admitting a reduction of type IV in the sense of \cite{Umi25}. 
Furthermore, in Section~3  we define essential reductions of the tame automorphisms   
and study the uniqueness of these reductions.  In particular, it is proved that every tame automorphism admits 
only one type of the essential reductions. A set of defining relations of the tame automorphism group $TA_3(F)$ is given in Section~4.

The results of this paper were first published in \cite{UUU}.

\section{Reductions of tame automorphisms}

\hspace*{\parindent}

Let  $F$ be an arbitrary field of characteristic $0$, and let  
$A=F[x_1,x_2,x_3]$ be the polynomial algebra in the variables  $x_1,x_2,x_3$ over  $F$. 
As in \cite{Umi24}, we will identify the algebra  $A$ with the corresponding subspace 
of the free Poisson algebra   $P=P\langle x_1,x_2,x_3\rangle$. 
The highest homogeneous part  $\overline{f}$ and the degree  $\deg f$ can be defined in an ordinary way. 
Note that  
\bes 
\overline{fg}=\overline{f}
\overline{g},\ \deg(fg)=\deg f+\deg g,\ \deg[f,g]\leq \deg f+\deg g.
\ees

If  $f_1,f_2,\ldots ,f_k\in A$, then we denote by  
$\langle f_1,f_2,\ldots,f_k\rangle$  the subalgebra of  $A$ generated by these elements.  

We will often use the terminology and results of  \cite{Umi25}. Therefore, let us make  
 an agreement that Corollary 2A, Lemma 4A, 
 Proposition  5A, and Theorem  7A mean Corollary~2, Lemma 4, Proposition 5, and Theorem  7 of \cite{Umi25},
 respectively.

Let us remind now the more necessary definitions and statements from \cite{Umi24,Umi25}. 

\smallskip

\begin{lm}\label{l1.1}{\bf \cite{Umi24}}
Let  $f,g,h\in A$. Then the following statements are true: 

$(1)$ $[f,g]=0$ if and only if  $f$ and $g$ are algebraically dependent. 

$(2)$ Suppose that   $f,g,h \notin F$ and  $m=\deg[f,g]+\deg h$, $n=\deg[g,h]+\deg f$,
$k=\deg[h,f]+\deg g$.
Then  $m\leq \max(n,k)$. If $n\neq k$, then $m=\max(n,k)$.
\end{lm}

The following statements are well known (see  \cite{Cohn}):

\smallskip

(F1) {\em If $a$ and $b$  are homogeneous algebraically dependent elements of the algebra  $A$, 
then there exists an element  $z\in A$ such that   
$a=\a z^n$, $b=\b z^m$ and $\a,\b\in F$. The subalgebra 
$\langle a,b \rangle$ is single generated if and only if  $m|n$ or  $n|m$}. 

\smallskip

(F2) {\em Let  $f,g\in A$ be such that $\overline{f}$ and $\overline{g}$ are algebraically independent. If   
$h\in \langle f,g \rangle$, then  $\overline{h}\in \langle\overline{f},\overline{g}\rangle$.}

\medskip

A pair of elements  $f,g$ of the algebra  $A$ is called  
{\em reduced} if  $\overline{f}\notin
\langle\overline{g}\rangle $ and $\overline{g}\notin \langle\overline{f}\rangle$. A reduced pair 
of algebraically independent elements  $f,g\in A$ is called {\em $\ast$-reduced} 
if $\overline{f}$ and $\overline{g}$ are algebraically dependent.  

\smallskip

Consider a $\ast$-reduced pair of elements  $f,g$ of the algebra  $A$ and let  $n=\deg
f<m=\deg g$. Put  $p=\frac{n}{(n,m)}$, $s=\frac{m}{(n,m)}$, and 
\bes
N=N(f,g)=\frac{mn}{(m,n)}-m-n+\deg[f,g]=mp-m-n+\deg[f,g],
\ees
where  $(n,m)$ is the greatest common divisor of   $n$ and $m$. Note that  $(p,s)=1$, and since 
 $\overline{f}$ and $\overline{g}$ are algebraically dependent, there exists an element  
 $a\in A$ such that  $\overline{f}=\b a^p$ and 
$\overline{g}=\g a^s$. Sometimes we will call a  $\ast$-reduced pair of elements  $f,g$ also a 
 {\em $p$-reduced} pair. Assume that $G(x,y)\in
F[x,y]$. It was proved in  \cite{Umi24} that  if  $\deg_y(G(x,y))=pq+r$,
$0\leq r<p$, then  
\bee\label{f1.1}
\deg(G(f,g))\geq qN+mr,
\eee
and if  $\deg_x(G(x,y))=sq_1+r_1$, $0\leq r_1<s$, then  
\bee\label{f1.2}
\deg(G(f,g))\geq q_1N+nr_1.
\eee

\smallskip

\begin{co}\label{c1.1}{\bf \cite{Umi25}}
Assume that $G(x,y)\in F[x,y]$ and $h=G(f,g)$. Consider the following conditions: 
\begin{itemize}
\item[$(i)$] $\deg h<N(f,g)$;
\item[$(ii)$] $\deg_y(G(x,y))<p$;
\item[$(iii)$] $h=\sum_{i,j}\a_{ij}f^ig^j$, where $\a_{ij}\in F$
 and  $in+jm\leq\deg h$ for all  $i,j$;
\item[$(iv)$] $\overline{h}\in \langle\overline{f},\overline{g}\rangle$.
\end{itemize}
Then  $(i)\Rightarrow (ii)\Rightarrow (iii)\Rightarrow (iv)$.
\end{co}

\smallskip

\begin{lm}\label{l1.2}{\bf \cite{Umi25}}
There exists a polynomial  $w(x,y)\in F[x,y]$ of the type 
\bes
w(x,y)=y^p-\a x^s-\sum \a_{ij}x^iy^j, \ ni+mj<mp,
\ees
which satisfies the following conditions: 

$(1)$ $\deg w(f,g)<pm$;

$(2)$ $\overline{w(f,g)}\notin \langle\overline{f},\overline{g}\rangle$.
\end{lm}

A polynomial $w(x,y)$ satisfying the conditions of Lemma  \ref{l1.2} is called a  {\em derivative polynomial} of the  $\ast$-reduced pair $f,g$.

\begin{co}\label{c1.2}{\bf \cite{Umi25}}
If  $h\in \langle f,g\rangle \setminus F$ and 
$\deg h<n$, then $h=\l
w(f,g),\ 0\neq \l\in F,$ where $w(x,y)$ is a derivative polynomial of the pair $f,g$.
\end{co}

\begin{co}\label{c1.3}{\bf \cite{Umi25}}
If $w(x,y)$ is a derivative polynomial of the pair $f,g$, then  
\bes
\deg(\frac{\d w}{\d x}(f,g))=n(s-1), \\
\deg(\frac{\d w}{\d y}(f,g))=m(p-1).
\ees
\end{co}

\smallskip

Let $\te =(f_1,f_2,f_3)$ be an arbitrary automorphism of the algebra  $A$. The number  
\bes
\deg \te=\deg f_1+\deg f_2+\deg f_3
\ees
is called the {\em degree} of $\te$.

Recall that an {\em elementary transformation} of a triple  
$(f_1,f_2,f_3)$ is, by definition, a transformation that changes only one element $f_i$ to an element of the form $\a
f_i+g$, where  $0\neq \a\in F$ and $g\in \langle\{f_j|j\neq
i\}\rangle$. The notation  
\bes
\te \rightarrow \tau
\ees 
means that the triple $\tau$ is obtained from $\te$ by a single elementary transformation. 
An automorphism $\te$ is called {\em  elementarily reducible} or {\em admits an elementary reduction} if there exists
     $\tau\in GA_3(F)$ such that $\te\rightarrow \tau$ and $\deg\tau<\deg\te$. The element  $f_i$ of the automorphism $\te$ which was changed in $\tau$   
 to an element of less degree is called  
{\em reducible} and we will say also that  $f_i$ is {\em reduced in $\te$ by the automorphism $\tau$}.

\begin{df}\label{d1.1}{\bf \cite{Umi25}}
Let  $\te=(f_1,f_2,f_3)$ be an automorphism of $A$ such that $\deg f_1=2n$,
$\deg f_2=ns$, $s\geq 3$ is an odd number, $2n<\deg f_3\leq ns$, and 
$\overline{f}_3\notin \langle\overline{f}_1,\overline{f}_2\rangle$. Suppose that there exists  
 $0\neq \b \in F$ such that the elements $g_1=f_1$ and $g_2=f_2-\b f_3$ satisfy the conditions:
\begin{itemize}
\item[(i)] $g_1,g_2$ is a $2$-reduced pair and $\deg g_1=\deg f_1$, $\deg g_2=\deg f_2$;
\item[(ii)] the element  $f_3$ of the automorphism $(g_1,g_2,f_3)$ is reduced by an automorphism  
 $(g_1,g_2,g_3)$ with the condition  
$\deg[g_1,g_3]<\deg g_2+\deg[g_1,g_2]$.
\end{itemize}
Then we will say that $\te$ admits a reduction of type I and the automorphism $(g_1,g_2,g_3)$ 
will be called a reduction of $\te$ of type  I with an active element $f_3$.
\end{df}

\smallskip

\begin{df}\label{d1.2}{\bf \cite{Umi25}}
Let  $\te=(f_1,f_2,f_3)$ be an automorphism of $A$ such that $\overline{f}_1$ and $\overline{f}_3$ are linearly independent,  $\deg f_1=2n$,
 $\deg f_2=3n$, and $\frac{3n}{2}<\deg f_3 \leq 2n$. 
Suppose that there exist 
$\a,\b \in F$,  where $(\a,\b)\neq 0$, such that the elements  
$g_1=f_1-\a f_3$ and $g_2=f_2-\b f_3$ satisfy the conditions (i) and (ii)   of Definition  \ref{d1.1}. 
Then we will say that $\te$ admits a reduction of type II and the automorphism $(g_1,g_2,g_3)$ 
will be called a reduction of $\te$ of type  II with an active element $f_3$.
\end{df}

\smallskip

\begin{df}\label{d1.3}{\bf \cite{Umi25}}
Let $\te=(f_1,f_2,f_3)$ be an automorphism of $A$ such that $\deg f_1=2n$, and either $\deg f_2=3n$, $n<\deg f_3 \leq \frac{3n}{2}$, or  
$\frac{5n}{2}<\deg f_2\leq 3n$, $\deg f_3 =\frac{3n}{2}$. Suppose that there exist 
$\a,\b_1,\b_2 \in F$, where $(\a,\b_1,\b_2) \neq 0$, such that the elements  
$g_1=f_1-\a f_3$ and $g_2=f_2-\b_1 f_3-\b_2 f_3^2$ satisfy the conditions: 
\begin{itemize}
\item[(i)] $g_1,g_2$ is a $2$-reduced pair and $\deg g_1=2n$,
$\deg g_2=3n$;
\item[(ii)] the element $f_3$ of the automorphism $(g_1,g_2,f_3)$   is reduced by an automorphism  $(g_1,g_2,g_3)$ with the condition 
$\deg g_3<n+\deg[g_1,g_2]$.
\end{itemize}
Then we will say that $\te$ admits a reduction of type III and the automorphism $(g_1,g_2,g_3)$ 
will be called a reduction of $\te$ of type  III with an active element $f_3$. 
\end{df}

\smallskip
The next remark can be extracted from the proofs of Propositions 1A, 2A, and 3A. 
\begin{rem}\label{r1.1}
\begin{itemize}
\item[(i)] $\deg\,[g_1,g_2]\leq 2n$ in the conditions of Definition \ref{d1.1}; 
\item[(ii)] $\deg\,[g_1,g_2]\leq n$ in the conditions of Definition \ref{d1.2}; 
\item[(iii)] $\deg\,[g_1,g_2]\leq \frac{n}{2}$ in the conditions of Definition \ref{d1.3}.
\end{itemize}
\end{rem}

\begin{lm}\label{l1.3}
Let $\te=(f_1,f_2,f_3)$ satisfy the conditions of Definition \ref{d1.3}. 
Then the following statements are true:
\begin{itemize}
\item[(i)] if $a\in <f_2,f_3>$ and  $\deg\,a\leq \frac{5n}{2}$, then  $a\in <f_3>$; 
\item[(ii)] if $a\in <f_1,f_3>$ and $\deg\,a\leq \frac{7n}{2}$, then $\overline{a}\in <\overline{f_1},\overline{f_3}>$; 
\item[(iii)] if $a\in <f_1,f_2>$ and $\deg\,a<2n$, then $a= \delta\in F$. 
\end{itemize}
\end{lm}

We omit the proof of this lemma, since later in Lemma \ref{l1.4} we consider analogous statements for a more complicated case. 

The next definition is some extension of the definition of automorphisms admiting a reduction of type IV in the sense of \cite{Umi25}. 

\begin{df}\label{d1.4}
Let $\te=(f_1,f_2,f_3)$ be an automorphism of $A$ such that  
$\deg\,f_3\leq\frac{3n}{2}$ and $\deg\,\te \leq \frac{13n}{2}$. Suppose that there exist  $\a_1,\a_2,\b_1,\b_2,\b_3,\b_4 \in F$ such that the elements  
 $g_1=f_1-\a_1 f_3-\a_2 f_3^2$ and $g_2=f_2-\b_1 f_3-\b_2 f_3^2-\b_3 f_1 f_3-\b_4 f_3^3$ 
 satisfy the conditions: 
\begin{itemize}
\item[(i)] $g_1,g_2$ is a $2$-reduced pair and $\deg g_1=2n$,
$\deg g_2=3n$;
\item[(ii)] there exists an element $g_3$ of the form  
\bes
g_3= f_3-\g w(g_1,g_2), \,\, 0\neq \g \in F, 
\ees
where $w(x,y)$ is a derivative polynomial of the 2-reduced pair  $g_1,g_2$, such that 
\begin{itemize}
\item[(a)] $\deg g_3=\frac{3n}{2}$ and $\deg[g_1,g_3]<3n+\deg[g_1,g_2]$; 
\item[(b)] there exists $0\neq \mu \in F$ such that $\deg(g_2-\mu g_3^2)\leq 2n$. 
\end{itemize}
\end{itemize}
Then we will say that $\te$ admits a reduction of type IV and the automorphism $(g_1,g_2-\mu g_3^2,g_3)$ 
will be called a reduction of $\te$ of type  IV with an active element $f_3$. 
\end{df}

We will also use Definitions \ref{d1.1}--\ref{d1.4}, admitting a permutation of the components of $(f_1,f_2,f_3)$. 

It is difficult to find examples of automorphisms illustrating the definitions \ref{d1.1}--\ref{d1.4}. 
An example of an automorphism admits a reduction of type I was constructed in \cite[p. 204, Example 1]{Umi25}. 
At the moment we have no example of an automorphism admits a reduction of types II--IV and the corresponding question was also formulated in  \cite[p. 225, Problem 1]{Umi25}. 

\begin{pr}\label{p1.1}
Let $\te=(f_1,f_2,f_3)$ be an automorphism of $A$ satisfying the conditions of Definition \ref{d1.4}. Then the following statements are true:  
\begin{itemize}
\item[(1)] $\deg\,[g_1,g_2]\leq\frac{n}{2}$, $\deg\,f_3\geq n+\deg\,[g_1,g_2]$; 
\item[(2)] $2n\leq\deg\,f_1<\frac{5n}{2}$, $\frac{5n}{2}+\deg\,[g_1,g_2]\leq \deg\,f_2<\frac{7n}{2}$;
\item[(3)] $\deg\,[f_1,f_3]=\deg\,[g_1,f_3]>3n$, $\deg\,[f_2,f_3]>3n$; 
\item[(4)] $\deg\,f_2+\deg\,f_3>4n$; 
\item[(5)] if $(\b_3,\b_4)=0$, then either $\deg\,f_2=3n$ or $\deg\,f_3=\frac{3n}{2}$; 
\item[(6)] the coefficients $\a_1,\a_2,\b_1,\b_2,\b_3,\b_4$ are uniquely defined; 
\item[(7)] if $(\a_1,\a_2,\b_1,\b_2,\b_3,\b_4)\neq 0$, then  $\deg\,[f_1,f_2]>3n$. 
\end{itemize}
\end{pr}
\Proof
Since  $\deg\,g_3, \deg\,f_3\leq \frac{3n}{2}$, and $\g\neq 0$, Definition \ref{d1.4} 
gives $\deg\,w(g_1,g_2)\leq \frac{3n}{2}$. By the definition of a derivative polynomial, applying Corollary \ref{c1.1} we get  
\bes
\frac{3n}{2}\geq \deg\,w(g_1,g_2)\geq N(g_1,g_2)=n+\deg\,[g_1,g_2]. 
\ees
Consequently, $\deg\,[g_1,g_2]\leq \frac{n}{2}$. 
According to Corollary \ref{c1.3} we have  $\deg\,\frac{\partial w}{\partial y}(g_1,g_2)=3n$. 
By Definition \ref{d1.4}, we have also   
\bes
\deg\,[g_1,g_3]<3n+\deg\,[g_1,g_2]=\deg\,[g_1,g_2]\frac{\partial w}{\partial y}(g_1,g_2).  
\ees
Since   
\bes
[f_1,f_3]=[g_1,f_3]=[g_1,g_3]+[g_1,g_2]\frac{\partial w}{\partial y}(g_1,g_2),  
\ees
comparing the degrees of elements here we find that   
\bes
\deg\,[f_1,f_3]=\deg\,[g_1,f_3]=3n+\deg\,[g_1,g_2].
\ees
Since $\deg\,g_1=2n$, it follows also that $\deg\,f_3\geq n+\deg\,[g_1,g_2]$. Consequently,  $\deg\,f_3<\deg\,g_1=2n<\deg\,f_3^2$, which gives $\deg\,f_1\geq 2n$. 

Obviously,   
\bes
\deg\,[g_1,g_2]+\deg\,f_3< \deg\,[g_1,f_3]+\deg\,g_2. 
\ees
Then by Lemma \ref{l1.1}(2) we obtain  
\bes
\deg\,[g_2,f_3]+\deg\,g_1=\deg\,[g_1,f_3]+\deg\,g_2,  
\ees
i.e. 
\bes
\deg\,[g_2,f_3]=4n+\deg\,[g_1,g_2].   
\ees
Next we have  
\bes
\deg\,[g_2,f_3]<4n+2\deg\,[g_1,g_2]\leq \deg\,[f_1,f_3]f_3. 
\ees
Since   
\bes
[f_2,f_3]=[g_2,f_3]+\b_3 [f_1,f_3]f_3,  
\ees
comparing the degrees of elements we get  
\bes
\deg\,[f_2,f_3]\geq 4n+\deg\,[g_1,g_2].
\ees
Then 
\bes
\deg\,f_2\geq 4n+\deg\,[g_1,g_2]-\deg\,f_3\geq \frac{5n}{2}+\deg\,[g_1,g_2]. 
\ees

Assume that  $(\b_3,\b_4)\neq 0$. Then it is easy to check that $\deg\,(\b_3f_1+\b_4f_3^2)f_3>3n$. 
Since $\deg\,f_3<\deg\,f_3^2\leq \deg\,g_2=3n$, we have  $\deg\,f_2>3n$. Now if $(\b_3,\b_4)= 0$, then we have  
 $\deg\,f_2\leq 3n$. And the inequality  $\deg\,f_2< 3n$ is possible only if $\deg\,f_3=\frac{3n}{2}$. 
 It remains to note that the inequality
    $\deg\,f_2+\deg\,f_3>4n$ is fulfilled in both cases. 
    Since   $\deg\,\te\leq \frac{13n}{2}$, this gives also $\deg\,f_1<\frac{5n}{2}$. 
    Note that the inequality $\deg\,f_1+\deg\,f_3>3n$ implies that $\deg\,f_2<\frac{7n}{2}$.

We have  
\bes
f_1=g_1+\a_1f_3+\a_2f_3^2, \\
f_2=g_2+\b_1f_3+\b_2f_3^2+\b_3f_1f_3+\b_4f_3^3\\
= g_2+\b_1f_3+(\b_2+\a_1\b_3)f_3^2+\b_3g_1f_3+(\b_4+\a_2\b_3)f_3^3. 
\ees
Consequently,  
\bes
[f_1,f_2]=[g_1,g_2]+\b_1[g_1,f_3]+\a_1[f_3,g_2]+(2\b_2+\a_1\b_3)[g_1,f_3]f_3\\
+ \b_3g_1[g_1,f_3]+2\a_2f_3[f_3,g_2]+(3\b_4+\a_2\b_3)[g_1,f_3]f_3^2.
\ees
Note that   
\bes
\deg\,[g_1,g_2]<3n<\deg\,[g_1,f_3]<\deg\,[f_3,g_2]<\deg\,[g_1,f_3]f_3\\
< \deg\,g_1[g_1,f_3]<\deg\,[f_3,g_2]f_3<\deg\,[g_1,f_3]f_3^2.  
\ees
Consequently, $3\b_4+\a_2\b_3=0$ or 
\bes
\overline{[f_1,f_2]}=(3\b_4+\a_2\b_3)\overline{[g_1,f_3]f_3^2}=(3\b_4+\a_2\b_3)\overline{[f_1,f_3]f_3^2}, 
\ees
i.e. $3\b_4+\a_2\b_3$ is uniquely defined. Now as in the proofs of Propositions 1A, 2A, and 3A we can easily deduce the statement (6) of Proposition  \ref{p1.1}. If $\deg\,[f_1,f_2]\leq 3n$, then obviously 
\bes
\b_1=\a_1=2\b_2+\a_1\b_3=\b_3=2\a_2=3\b_4+\a_2\b_3=0, 
\ees
i.e. $(\a_1,\a_2,\b_1,\b_2,\b_3,\b_4)=0$. $\Box$

\begin{lm}\label{l1.4}
Let  $\te=(f_1,f_2,f_3)$ satisfy the conditions of Definition \ref{d1.4}. 
Then the following statements are true:
\begin{itemize}
\item[(i)] if $a\in <f_2,f_3>$ and $\deg\,a\leq \frac{5n}{2}$, then $a\in <f_3>$; 
\item[(ii)] if $a\in <f_1,f_3>$ and $\deg\,a\leq \frac{7n}{2}$, then  $\overline{a}\in <\overline{g_1},\overline{f_3}>$;  
\item[(iii)] if  $a\in <f_1,f_2>$ and  $\deg\,a<2n$, then we have $a= \delta\in F$ if $(\a_1,\a_2,\b_1,\b_2,\b_3,\b_4)\neq 0$, while 
 $a= \delta w(g_1,g_2)+\lambda$ otherwise, where $w(x,y)$ is a derivative polynomial of the 2-reduced pair $g_1,g_2$. 
\end{itemize}
\end{lm}
\Proof
Assume that $a\in <f_2,f_3>$. The subalgebra $<f_2,f_3>$ does not change if we replace $f_2$ by $g_2+\b_3g_1f_3$. So, we can assume that 
$f_2=g_2+\b_3g_1f_3$.  If $\overline{f_2}$ and $\overline{f_3}$  
 are algebraically independent, then by  (F2) we have $\overline{a}\in <\overline{f_2},\overline{f_3}>$. 
By Proposition \ref{p1.1}, we have $\deg\,f_2>\frac{5n}{2}$. Consequently, $a\in <f_3>$. Suppose that $\overline{f_2}$ and $\overline{f_3}$
 are algebraically dependent. If  $\b_3\neq 0$, then $f_2,f_3$ is a  
 $\ast$-reduced pair. If $\b_3= 0$, then either $f_2,f_3$ is a  $\ast$-reduced pair, 
 or there exists $\mu \in F$ such that $\overline{g_2}=\mu \overline{f_3^2}$. 
 In the last case we can assume that $f_2=g_2-\mu f_3^2$. Then again $f_2,f_3$ is a $\ast$-reduced pair. Since  
\bes
N(f_2,f_3)> \deg\,[f_2,f_3]> 4n>\frac{5n}{2}\geq \deg\,a, 
\ees 
the inequality $\deg\,a\geq N(f_2,f_3)$ leads to a contradiction.  Consequently,  \\
$\deg\,a< N(f_2,f_3)$.  
Since $\deg\,f_2>\frac{5n}{2}\geq \deg\,a$, by Corollary \ref{c1.1}  we get again $a\in <f_3>$.

Suppose that $a\in <f_1,f_3>$. As above we can assume that $f_1=g_1$. 
If $\overline{g_1}$ and $\overline{f_3}$     
 are algebraically dependent, then $g_1,f_3$ is a $p$-reduced pair. 
Consider the case when  
\bes
\frac{7n}{2}\geq \deg\,a\geq N(g_1,f_3)= (p-1)\deg\,g_1 -\deg\,f_3+\deg[g_1,f_3].
\ees
If  $p\geq 3$, then with regard to Proposition \ref{p1.1} we have $N(g_1,f_3)>5n$. If $p=2$, then this is possible only if  
\bes
\deg\,g_1=2n=3k, \,\, \deg\,f_3=2k\leq \frac{3n}{2}.
\ees
Consequently, repeately applying Proposition \ref{p1.1} gives  
\bes
N(g_1,f_3)= k+\deg[g_1,f_3]>\frac{2n}{3}+3n=\frac{11n}{3}>\frac{7n}{2}.
\ees
Thus $\deg\,a < N(g_1,f_3)$. Then Corollary \ref{c1.1} gives $\overline{a}\in <\overline{g_1},\overline{f_3}>$. 

Now suppose that $a\in <f_1,f_2>$ and $\deg\,a<2n$. If  $(\a_1,\a_2,\b_1,\b_2,\b_3,\b_4)\neq 0$, then
 by Proposition \ref{p1.1} we have $\deg\,[f_1,f_2]>3n$. By Corollary \ref{c1.1}, from here we deduce that $\overline{a}\in <\overline{f_1},\overline{f_2}>$. Since  $\deg\,f_1, \deg\,f_2\geq 2n$, it follows that  
 $a=\delta \in F$. If $(\a_1,\a_2,\b_1,\b_2,\b_3,\b_4)= 0$, then Corollary \ref{c1.2} gives the 
 statement (iii) of the lemma. $\Box$

\smallskip

It is easy to deduce from Proposition \ref{p1.1} that  if $(\a_2,\b_3,\b_4)= 0$, then  
 Definition \ref{d1.4} gives exactly the automorphisms admitting a reduction of type IV in the sense of \cite{Umi25}. 
\begin{theor}\label{t1.1}{\bf \cite{Umi25}}
Let $\te=(f_1,f_2,f_3)$ be a tame automorphism of the polynomial algebra    
$A=F[x_1,x_2,x_3]$ over a field $F$ of characteristic $0$.  If  $\deg\te>3$, then $\te$ admits 
either an elementary reduction or a reduction of one of the types I--IV.
\end{theor}

\smallskip
Further we need the following proposition.  
\begin{pr}\label{p1.2} 
Let $(f_1,f_2,f_3)$ be a tame automorphism of the algebra $A$   
 satisfying the following conditions:
\begin{itemize}
\item[(i)] $f_1,f_2$ is a $2$-reduced pair and $\deg f_1=2n$, $\deg f_2=sn$, where $s\geq 3$ is an odd number; 
\item[(ii)]  $\deg f_3=m<sn$ and  $f_3$ is an irreducible element of the automorphism $(f_1,f_2,f_3)$.
\end{itemize}
Then one of the following statements is true: 
\begin{itemize}
\item[(1)] $m<n(s-2)+\deg[f_1,f_2]$;
\item[(2)] $(f_1,f_2,f_3)$ admits a reduction of type  IV with the active element $f_3$, and the coefficients 
 $\a_1,\a_2,\b_1,\b_2,\b_3,\b_4$ of Definition \ref{d1.4} are equal to $0$;
\item[(3)] $\deg[f_1,f_3]<3n+\deg[f_1,f_2]$ and there exists $\mu\in F$
 such that $\deg(f_2-\mu f_3^2)\leq 2n$.
\end{itemize}
\end{pr}
\Proof
The conditions of this proposition coincide with the conditions of Propositions 4A and  5A  
if we take into account the fact that the tame automorphisms of the algebra $A$ are already 
simple  and every elementarily reducible element is simple reducible in the sense of \cite{Umi25}. 
It is easy to check that the proofs of Propositions 4A and 5A give also the proof of Proposition \ref{p1.2}. 
$\Box$

\smallskip

\section{On the uniqueness of reductions}

\hspace*{\parindent}

Every tame automorphism  $\te$ has a sequence of elementary transformations  

\bee\label{f2.1}
(x_1,x_2,x_3)= \te_0\rightarrow
\te_1\rightarrow \ldots \rightarrow
\te_{k-1}\rightarrow \te_k=\te.
\eee

Put  $d=\max\{\deg\,\te_i \,|\, 0\leq i\leq k\}$. The number  $d$ will be called 
the {\em width} of the sequence  (\ref{f2.1}). The minimal width of all sequences 
of the type (\ref{f2.1}) for $\te$ will be  called the {\em width} of the automorphism $\te$. 

Let $\te \in TA_3(F)$, $\deg\,\te>3$. An automorphism $\phi$ will be called 
an {\em essential reduction} of $\te$ if there exists a 
sequence of elementary transformations  

\bee\label{f2.2}
\phi \rightarrow \te_0\rightarrow
\te_1\rightarrow \ldots \rightarrow
\te_{t-1}\rightarrow \te_t=\te 
\eee
such that  $d(\phi)<d(\theta)$ and  $\deg\,(\te_i)\leq d(\theta)$, where $0\leq i\leq t$. 

The minimal number $t$ for which there exists an essential reduction 
 $\phi$ of $\theta$ of  the type (\ref{f2.2}) will be called   the {\em height} 
 of $\theta$, and the corresponding sequence  
 (\ref{f2.2}) will be called a {\em minimal essential reduction} 
of  $\theta$. If  $t=0$, then we will say that $\te$ admits an {\em essential elementary reduction}. 

Later we will see that the elementary reductions of an automorphism admitting a reduction of type  III or IV  are not essential.  

\smallskip

If $\deg\,\te=3$, then we put $t=\infty$. Introduce a lexicographic order  
 on the set of all pairs $(d,t)$, where $d$ is the width and $t$ is the height of some tame automorphism, by putting $(d_1,t_1)<(d_2,t_2)$ if either 
 $d_1<d_2$ or $d_1=d_2,\,t_1<t_2$. 
 
\begin{lm}\label{l2.1} 
If (\ref{f2.2}) is a minimal essential reduction of $\theta$, then 
  $\deg\,\te_0=d(\te)$ 
and 
\bes
(d(\phi),t(\phi))<(d(\te_0),t(\te_0))<(d(\te_1),t(\te_1))<\ldots <(d(\te_t),t(\te_t))=(d(\te),t(\te)). 
\ees
\end{lm}
\Proof
If $\deg\,\te_0<d(\te)$, then (\ref{f2.2}) gives $d(\te_0)<d(\te)$. In this case instead of $\phi$ we can take $\te_0$, which contradicts   
the  minimality of (\ref{f2.2}). In addition, by the minimality of (\ref{f2.2}), $\phi$ is also a minimal 
essential reduction of each $\te_i,\,\,0\leq i\leq k$. Consequently,  
\bes
\deg\,\te_0= d(\te_i),\,\, 0\leq i\leq k, \,\,\, t(\te_{i+1})=t(\te_i)+1,\,\,0\leq i\leq t-1, 
\ees
which gives the statement of the lemma.  
$\Box$

Let $d$ be the width and $t$ be the height of $\theta$. The sequence 
  (\ref{f2.1}) will be called a {\em minimal representation} of the automorphism $\theta$ if  
\bes
\te_{k-t-1} \rightarrow \te_{k-t}\rightarrow
\te_{k-t+1} \rightarrow \ldots \rightarrow \te_k=\te
\ees
 is a minimal essential reduction of $\theta$ and  $\deg\,\te_i<d(\te),\,\,0\leq i\leq k-t-2$. 

So, every minimal essential reduction of a tame automorphism can be continued to a minimal representation. 
\smallskip

\begin{pr}\label{p2.1}
Let $\te\in TA_3(F)$, $\deg \te>3$, let $d$  be the width and  $t$ be the height 
 of $\te$, and let (\ref{f2.2}) be a minimal essential reduction of $\te$. Then the following statements hold:  
\begin{itemize}
\item[(a)] If $\te$ admits a reduction of type I, then $d=\deg\,\te$, $t=1$, and in the conditions of Definition \ref{d1.1} we have  
\bee\label{f2.3} 
\te_{t-1}=(f_1, \eta g_2+\kappa g_1 + \nu, f_3), \,\,\eta \neq 0.  
\eee

\item[(b)] 
If  $\te$ admits a reduction of type II, then $d=\deg\,\te$ and in the conditions 
of Definition \ref{d1.2} the automorphism $\te_{t-1}$ and the height $t$ will be calculated in the following way:  
\begin{itemize}
\item[(1)]  if $\a \b \neq 0$, then $t=2$ and  $\te_{t-1}$ can be written down simultaneously in the form (\ref{f2.3}) and in the form    
\bee\label{f2.4} 
\te_{t-1}=(\rho g_1 + \sigma, f_2, f_3), \,\,\rho \neq 0;  
\eee
\item[(2)]  if $\a=0$, then  $t=1$ and  $\te_{t-1}$ has the form (\ref{f2.3}); 
\item[(3)]  if $\b=0$, then $t=1$ and $\te_{t-1}$ has the form (\ref{f2.4}). 
\end{itemize}
\item[(c)]
If $\te$ admits a reduction of type III, then in the conditions of Definition \ref{d1.3} 
we have $d=5n+\deg\,f_3$  and  the automorphism  $\te_{t-1}$ and the height $t$ will be calculated in the following way:  
\begin{itemize}
\item[(1)]  if $\a\neq 0$ and $(\b_1,\b_2) \neq 0$, then $t=2$ and $\te_{t-1}$ can be written down simultaneously 
in the form (\ref{f2.3}) and (\ref{f2.4}); 
\item[(2)]  if $\a = 0$, then $t=1$ and $\te_{t-1}$ has the form (\ref{f2.3}); 
\item[(3)] if $(\b_1,\b_2) = 0$, then $t=1$ and $\te_{t-1}$ has the form (\ref{f2.4}). 
\end{itemize}
\item[(d)]
If $\te$ admits a reduction of type IV, then in the conditions of Definition \ref{d1.4} we have $d=\frac{13n}{2}$ and  
the automorphism  $\te_{t-1}$  and the height $t$ will be calculated in the following way:   
\begin{itemize}
\item[(1)]  if $(\a_1,\a_2)\neq 0$ and $(\b_1,\b_2,\b_3,\b_4)\neq 0$, then $t=3$ and $\te_{t-1}$ can be written down simultaneously 
in the form (\ref{f2.3}) and (\ref{f2.4}); 
\item[(2)]  if $(\a_1,\a_2)=0$ and $(\b_1,\b_2,\b_3,\b_4)\neq 0$, then $t=2$ and $\te_{t-1}$ has the form (\ref{f2.3}); 
\item[(3)] if $(\a_1,\a_2)\neq 0$ and $(\b_1,\b_2,\b_3,\b_4)=0$, then $t=2$ and  $\te_{t-1}$ has the form  (\ref{f2.4});  
\item[(4)] if $(\a_1,\a_2)=0$ and $(\b_1,\b_2,\b_3,\b_4)=0$, then $t=1$ and 
 $\te_{t-1}=(g_1,g_2,\rho g_3+\sigma)$.  
\end{itemize}
\item[(e)]
In the other cases $\te$ admits an essential elementary reduction, i.e.  
 $d=\deg\,\te$  and  $t=0$.
\end{itemize}
\end{pr}

\Proof 
Let  $\te\in TA_3(F)$ and let  $d$ be the width and  $t$ be the height of $\te$. 
Denote by $(d_1,t_1)$ the value of $(d,t)$ indicated in Proposition \ref{p2.1}.  We first show that  
\bee\label{f2.5} 
  (d,t)\leq (d_1,t_1) 
\eee
by induction on $\deg\,\te$. So, we can assume that the inequality (\ref{f2.5}) is true for automorphisms of less degree.  

We next proceed with the case study. The cases when $\te$ satisfies one of the conditions 
(a), (b), (c), and (d) are similar. Therefore, we give a proof in the case (d). 
  In this case  $\te$ admits a reduction of type 
 IV, and we adopt all the conditions and notation of Definition \ref{d1.4}. 
 Consider the sequence of elementary transformations  
\bee\label{f2.6} 
  \phi=(g_1,g_2-\mu g_3^2,g_3)\longrightarrow (g_1,g_2,g_3)\longrightarrow (g_1,g_2,f_3)
  \longrightarrow (f_1,g_2,f_3)\longrightarrow (f_1,f_2,f_3). 
\eee
Since $\deg\,\phi<\deg\,\te$, it follows that the inequality (\ref{f2.5}) is valid for $\phi$.  

We will show that  
\bee\label{f2.7} 
d(\phi)< \deg\,\te.
\eee

First of all we will give a standard argument that deduces 
 (\ref{f2.5}) from (\ref{f2.7}). In the sequence (\ref{f2.6})
 the automorphism $(g_1,g_2,g_3)$ has the maximal degree and $\deg\,(g_1,g_2,g_3)=\frac{13n}{2}=d_1$. 
 If (\ref{f2.7}) is fulfilled, then we get $d\leq d_1$.  If $d<d_1$, then (\ref{f2.5})  
 is fulfilled. Suppose that $d=d_1$.  Then the sequence (\ref{f2.6}) gives that $\phi$ is an essential reduction of $\te$. Consequently, $t\leq 3$. 
 We have $f_1=g_1$ if $(\a_1,\a_2)=0$, and  $f_2=g_2$ if $(\b_1,\b_2,\b_3,\b_4)=0$. In these cases, excluding $(f_1,g_2,f_3)$ from 
 (\ref{f2.6}), we get $t\leq 2$. If $(\a_1,\a_2,\b_1,\b_2,\b_3,\b_4)=0$, then excluding also $(g_1,g_2,f_3)$ from (\ref{f2.6}), we get $t\leq 1$. Thus, in all the cases $t\leq t_1$, and (\ref{f2.5}) is true. 

If $\phi$ satisfies one of the conditions (a), (b), and (e) of  the  
 proposition, then by (\ref{f2.5}) we have $d(\phi)\leq \deg\,(\phi)$, which gives (\ref{f2.7}). 

We have  
\bee\label{f2.8} 
\deg\,g_1=2n,\,\,\deg\,(g_2-\mu g_3^2)\leq 2n,\,\, \deg\,g_3=\frac{3n}{2}. 
\eee

Suppose that $\phi$ admits a reduction of type III or IV. According to 
  (\ref{f2.8}), the element $g_1$ has the highest degree among the components of $\phi$. 
  Therefore the  active element of the reduction is $g_2-\mu g_3^2$ or  $g_3$. 
  If $g_3$ is the active element of the reduction, then there exists $m$ such that  
\bes
m<\deg\,g_3\leq \frac{3m}{2},\,\,\,2m\leq \deg\,(g_2-\mu g_3^2)<\frac{5m}{2},\,\,\,\frac{5m}{2}<\deg\,g_1<\frac{7m}{2}.
\ees
Taking account of (\ref{f2.8}), from here we deduce the inequalities  
\bes
2m\leq 2n,\,\,\,\frac{3n}{2}\leq \frac{3m}{2},\,\,\,\frac{5m}{2}<\deg\,g_1=2m,
\ees
which are in a contradiction.  

If $g_2-\mu g_3^2$ is the active element of the reduction of $\phi$, then 
 there exists $m$ such that  
\bes
m<\deg\,(g_2-\mu g_3^2)\leq \frac{3m}{2},\,\,\,2m\leq \deg\,g_3<\frac{5m}{2},\,\,\,\frac{5m}{2}<\deg\,g_1<\frac{7m}{2}.
\ees
Consequently, $2m\leq \frac{3n}{2}$, i.e. $m\leq \frac{3n}{4}$. 
By (\ref{f2.5}), we have $d(\phi)\leq \frac{13m}{2}\leq \frac{39n}{8}<5n$. From Proposition \ref{p1.1} we have $\deg\,\te>6n$, i.e. 
 the inequality  (\ref{f2.7}) is fulfilled.  

\medskip 

Now, consider the case (e). If $\te$ does not admit  reductions of types I--IV, 
then according to Theorem \ref{t1.1} $\te$ admits an elementary reduction. Let $\phi$ be 
an elementary reduction of $\te$. If 
 $\phi$ satisfies one of the conditions (a), (b), and (e), then by (\ref{f2.5}) we have 
 $d(\phi)\leq \deg\,\phi<\deg\,\te$, i.e. the inequality (\ref{f2.7}) 
 is fulfilled.  

Assume that $\phi$ admits a reduction of type IV. Temporarily  we assume 
that $\phi=(f_1,f_2,f_3)$ and that $(f_1,f_2,f_3)$ satisfies the conditions of 
Definition \ref{d1.4}. According to (\ref{f2.5}), we have $d(\phi)\leq \frac{13n}{2}$. 
Put also $\te=(q_1,q_2,q_3)$.  
 If $\deg\,\te>\frac{13n}{2}$, then (\ref{f2.7}) 
 is fulfilled. Therefore we may assume that  
\bee\label{f2.9} 
\deg\,\te \leq \frac{13n}{2}. 
\eee
We first consider the case in which   
\bee\label{f2.10} 
q_1=\rho f_1+a,\,\,q_2=f_2,\,\,q_3=f_3,\,\,\rho\neq 0,\,\,a\in<f_2,f_3>. 
\eee
Changing $\rho g_1$ to  $g_1$, we may assume that $\rho=1$.
Proposition \ref{p1.1} gives $\deg\,f_2+\deg\,f_3>4n$ and $\deg\,f_1<\frac{5n}{2}$. Then from (\ref{f2.9}) we obtain   
 $\deg\,q_1<\frac{5n}{2}$. Consequently, $\deg\,a<\frac{5n}{2}$. By Lemma \ref{l1.4} we have $a\in <f_3>$ and  
\bes
a=\a_0'+\a_1' f_3+\a_2' f_3^2.
\ees
Consequently,  
\bes
q_1=g_1+\a_0'+(\a_1+\a_1')f_3+(\a_2+\a_2')f_3^2.
\ees
Changing $g_1+\a_0'$ to $g_1$, we may assume that $\a_0'=0$. Then  
 $\te$ admits a reduction of type IV, which contradicts the condition (e).  
 
 Let  
\bee\label{f2.11} 
q_1=f_1,\,\,q_2=\eta f_2+a,\,\,q_3=f_3,\,\,\eta\neq 0,\,\,a\in<f_1,f_3>. 
\eee
As above we can take $\eta=1$. By Proposition \ref{p1.1} we have  
 $\deg\,f_1+\deg\,f_3>3n$ and $\frac{5n}{2}<\deg\,f_2<\frac{7n}{2}$. 
Applying (\ref{f2.9}) we find that $\deg\,q_2<\frac{7n}{2}$, i.e. $\deg\,a<\frac{7n}{2}$. 
By Lemma \ref{l1.4} we obtain $\overline{a}\in <\overline{g_1},\overline{f_3}>$. 
Now it is easy to deduce that  
\bes
a=\g_0+\g_1g_1+\b_1'f_3+\b_2'f_3^2+\b_3'f_1f_3+\b_4'f_3^3.
\ees 
Thus,  
\bes
q_2=g_2+\g_0+\g_1g_1+(\b_1+\b_1')f_3+(\b_2+\b_2')f_3^2+(\b_3+\b_3')f_1f_3+(\b_4+\b_4')f_3^3.
\ees 
Changing  $g_2+\g_0+\g_1 g_1$ to $g_2$, we may assume that  
 $\g_0=\g_1=0$. Then  $\te$  admits a reduction of type IV. 
 
Now we consider the case in which    
\bee\label{f2.12} 
q_1=f_1,\,\,q_2=f_2,\,\,q_3=\rho f_3+a,\,\,\rho\neq 0,\,\,a\in<f_1,f_2>. 
\eee
Changing $\rho g_3$ to $g_3$, we may assume that $\rho=1$. 
 Proposition \ref{p1.1} gives  $\deg\,f_1+\deg\,f_2>\frac{9n}{2}$. By (\ref{f2.9}), from here  
 we get $\deg\,q_3<2n$, i.e. $\deg\,a<2n$. If $(\a_1,\a_2,\b_1,\b_2,\b_3,\b_4)\neq 0$, 
 then according to Lemma \ref{l1.4}  we obtain $a=\delta \in F$, i.e. $\te$ admits a reduction 
 of type IV. If $(\a_1,\a_2,\b_1,\b_2,\b_3,\b_4)=0$, then $f_1=g_1$, $f_2=g_2$, and by 
 Lemma \ref{l1.4} we have  
\bes
a=\delta w(g_1,g_2)+\sigma.
\ees
Then  
\bes
q_3=g_3+(\gamma+\delta)w(g_1,g_2)+\sigma. 
\ees
Changing $g_3+\sigma$ to $g_3$, we can take $\sigma=0$. If 
 $\gamma+\delta \neq 0$, then $\te$ admits a reduction of type IV. 
 If $\gamma+\delta=0$, then $\te=(g_1,g_2,g_3)$ and instead of $\phi$ 
 we can take  
\bes
\phi'=(g_1,g_2-\mu q_3^2,q_3). 
\ees
Recall that, considering the case when $\te$
 admits a reduction of type IV, we have simultaneously proved that  
\bes
d(\phi')<\frac{13n}{2}.
\ees
Since $\deg\,(\te)=\frac{13n}{2}$, the inequality (\ref{f2.7}) is also fulfilled.  

 If $\phi$ admits a reduction of type III, then as above we can assume 
 that $\phi=(f_1,f_2,f_3)$ and that $(f_1,f_2,f_3)$ satisfies the conditions of Definition \ref{d1.3}. 
 This case can be settled by analogy to the case when $\phi$ admits a reduction of type IV. 
  We only note that instead of (\ref{f2.9}) in this case we have a stronger inequality  
\bes
\deg\,\te \leq 5n+\deg\,f_3. 
\ees

The inequality (\ref{f2.5}) is thus proved if $\te$ satisfies the conditions (e). 
 Therefore we assume that the inequality (\ref{f2.5}) is proved.  

\medskip

Now we begin a full proof of Proposition \ref{p2.1}. 
Assume that the statement of the proposition is not true. Let $\te$ be an automorphism with minimal $(d,t)$  
which does not satisfy the statement of the proposition.  
Note that if $\te$ satisfies the condition (e), then we have nothing to prove.   

Restrict ourselves to the case when $\te$ admits a reduction of type IV.  
Let  $\te$ satisfy the conditions of Definition \ref{d1.4}. Put $\te_{t-1}=
(q_1,q_2,q_3)$. Since  
\bee\label{f2.13} 
(d(\te_{t-1}),t(\te_{t-1}))<(d,t)  
\eee
by Lemma \ref{l2.1}, it follows that  Proposition 
 \ref{p2.1} is valid for $\te_{t-1}$. By  (\ref{f2.5}) we have $d=d(\te)\leq \frac{13n}{2}$. Then  
\bee\label{f2.14} 
d(\te_{t-1})\leq \frac{13n}{2}.   
\eee
Suppose that $(q_1,q_2,q_3)$ has the form (\ref{f2.10}). Repeating the same arguments  
which were given after (\ref{f2.10}), we can assume that  
\bes
q_1=g_1+(\a_1+\a_1')f_3+(\a_2+\a_2')f_3^2, 
\ees
i.e. $\te_{t-1}$ admits a reduction of type IV. 

Let $(\b_1,\b_2,\b_3,\b_4)\neq 0$. If $(\a_1+\a_1',\a_2+\a_2')\neq 0$, then according 
to the statement  (d) of Proposition \ref{p2.1} we have 
$d(\te_{t-1})=\frac{13n}{2}$ and $t(\te_{t-1})=3$. Since $\deg\,\te \leq \frac{13n}{2}$, 
from here we get $d=\frac{13n}{2}$ and $t=4$, which contradicts (\ref{f2.5}). 
Consequently, $(\a_1+\a_1',\a_2+\a_2')=0$. In this case Proposition \ref{p2.1} gives $d(\te_{t-1})=\frac{13n}{2}$ and $t(\te_{t-1})=2$. 
Hence, $d=\frac{13n}{2}$ and $t=3$. If $(\a_1,\a_2)\neq 0$, then Proposition \ref{p2.1} is valid.  
If $(\a_1,\a_2)= 0$, then the inequality (\ref{f2.5}) is not fulfilled. It means that $\te_{t-1}$ 
 is obtained from $\te$ by another elementary transformation!   

Now let  $(\b_1,\b_2,\b_3,\b_4)= 0$. As above we have $(\a_1+\a_1',\a_2+\a_2')=0$. 
Then $d(\te_{t-1})=\frac{13n}{2}$ and $t(\te_{t-1})=1$. 
Consequently, $d=\frac{13n}{2}$ and $t=2$. If $(\a_1,\a_2)\neq 0$, 
then Proposition \ref{p2.1} is valid, and if $(\a_1,\a_2)= 0$, 
then the inequality (\ref{f2.5}) is not fulfilled. 

The above discussion is standard and the other cases  
can be examined similarly. $\Box$

\smallskip

\begin{co}\label{c2.1}
Reductions of types I--IV are essential minimal reductions.  
\end{co}

\begin{co}\label{c2.2}
If $\te \in TA_3(F)$ admits a reduction of types I--IV, then the type of this reduction is uniquely defined.  
\end{co}
\Proof
Let $\te=(f_1,f_2,f_3)$ admit a reduction of types I--IV. By Proposition \ref{p2.1}, 
the active element $f_3$ of this reduction is uniquely characterized as a component  
of $\te$ which does not change at a   
 minimal reduction before appearing in an automorphism $(g_1,g_2,f_3)$ with the property $\deg\,[g_1,g_2]\leq \deg\,f_3$. 
 The last inequality is a generalization of the estimates of the degree of $[g_1,g_2]$ in Remark \ref{r1.1} and in Proposition \ref{p1.1}. 
 In addition, if $\deg\,f_1\leq \deg\,f_2$, then the roles of elements in Definitions \ref{d1.1}--\ref{d1.4}
 are uniquely defined.  

Now the reductions of types I, II, and III can be easily distinguished among themselves by the degree of $\deg\,f_3$.  

Assume that $\te$ admits a reduction of type IV and satisfies the conditions of Definition \ref{d1.4}. 
If $\a_2\neq 0$, then $\overline{f_1}=\a_2\overline{f_3}^2$. 
Consequently, $f_1$ is an elementarily reducible element of $\te$. 
Note that if $\te$ admits a reduction of types I--III, then $\te$ does not admit such an elementary reduction. 
If $\a_2=0$, then
$\deg\,f_1=2n$ and $\deg\,f_3\leq \frac{3n}{2}$. 
Consequently, $\te$ does not admit a reduction of type I or II. 
If $\te$ admits a reduction of type III, 
then the elements $g_1$ and $g_2$ are uniquely defined for both reductions. 
We have $g_3=f_3-\gamma w(g_1,g_2)$. 
Since $\deg\,[g_1,g_3]<3n+\deg\,[g_1,g_2]$ and $\deg\,[g_1,w(g_1,g_2)]=3n+\deg\,[g_1,g_2]$, 
it follows that $\g$ is also uniquely defined from the equality  
$\overline{[g_1,f_3]}=\g \overline{[g_1,w(g_1,g_2)]}$. Now if $\deg\,g_3<\frac{3n}{2}$, 
then $\te$ admits a reduction of type III, and if $\deg\,g_3=\frac{3n}{2}$, 
then $\te$ admits a reduction of type IV. So, 
$\te$ does not admit reductions of types III and IV simultaneously. $\Box$

As we can see from the proof of Corollary \ref{c2.2}, not only the type of a reduction 
of $\te$ but also the elements $g_1$ and $g_2$ are uniquely defined. 
 The element $g_3$ in Definitions \ref{d1.3} and  \ref{d1.4} is uniquely defined up to a 
 summand from the field $F$, since it is with this exactness that the element $w(g_1,g_2)$ is defined. 
 In Definition \ref{d1.1} instead of 
 $g_3$ we can take any element of the form $\rho g_3+H(g_1)$, where 
 $0\neq \rho \in F$ and $\deg\,H(g_1)\leq \deg\,g_3$. Furthermore, in Definition \ref{d1.2} instead of $g_3$ 
 we can take any element of the form $\rho g_3+\gamma$, where $0\neq \rho \in F$.

\smallskip

Proposition \ref{p2.1} and Corollary \ref{c2.2} immediately   give  the next stronger form of Theorem  \ref{t1.1}.  

\begin{theor}\label{t2.1}
Let $\te \in TA_3(F)$, $\deg\,\te>3$. Then $\te$ admits only one type of the  following reductions:  
an essential elementary reduction or a reduction of types I--IV.  
\end{theor}

\smallskip

Besides, from the proof of Proposition \ref{p2.1} we can extract the next corollary.  
\begin{co}\label{c2.3}
Let $\te \in TA_3(F)$, $\deg\,\te>3$. If $\te$ admits a reduction of types I--IV, 
then every elementary reduction of $\te$ again admits a reduction of the same type.  
\end{co}

In fact, the elementary reductions considered in this corollary are not essential. 

Note that if $\te$ admits a reduction of type I or II, then $\te$ does not admit any elementary reduction. 
If $\te$ admits a reduction of type III, then in the conditions of Definition \ref{d1.3}  
$\te$  admits an elementary reduction if and only if $\overline{f_2}=\varsigma \overline{f_3^2},\,\,\deg\,f_3=\frac{3n}{2}$, and $\deg\,f_2=3n$. 
Assume that $\te$ admits a reduction of type IV. If $(\a_2,\b_3,\b_4)=0$ in the conditions of Definition \ref{d1.4}, 
then $\te$ also admits an elementary reduction  if and only if $\overline{f_2}=\varsigma \overline{f_3^2},\,\,\deg\,f_3=\frac{3n}{2}$, and $\deg\,f_2=3n$.  
 If $(\a_2,\b_3,\b_4)\neq 0$, then $\te$ admits obvious elementary reductions.

\bigskip

\section{Defining relations of the group of tame automorphisms}

\hspace*{\parindent}

Let $A_n=F[x_1,x_2,\ldots,x_n]$ be the polynomial algebra over a field $F$ with the set of variables  $X=\{x_1,x_2,\ldots,x_n\}$. Then the group 
 $TA_n(F)$ is generated by all elementary automorphisms 
\bee\label{f3.1} 
\s(i,\a,f) = (x_1,\ldots,x_{i-1}, \a x_i+f, x_{i+1},\ldots,x_n),
\eee
where $0\neq\a\in F,\ f\in F[X\setminus \{x_i\}]$.

\smallskip 

Our aim in this section is to describe defining relations of the group $TA_n(F)$ 
with respect to the generators  (\ref{f3.1}).  It is easy to check that   
\bee\label{f3.2} 
\s(i,\a,f) \s(i,\b,g) = \s(i,\a \b, \b f+g). 
\eee

Note that from here we obtain trivial relations  
\bes
\s(i,1,0) = id, \ \ 1\leq i\leq n.   
\ees

\smallskip 

If $i\neq j$ and $f\in F[X\setminus \{x_i,x_j\}]$, then we have also 
\bee\label{f3.3} 
\s(i,\a,f)^{-1} \s(j,\b,g) \s(i,\a,f) = \s(j,\b,\s(i,\a,f)^{-1}(g)). 
\eee

Consequently, if $i\neq j$ and $f,g\in F[X\setminus \{x_i,x_j\}]$, then  the automorphisms $\s(i,\a,f)$,  
$\s(j,\b,g)$ commute.  

\smallskip

For every pair $k,s$, where $1\leq k\neq s\leq n$, we define a tame automorphism $(ks)$ by putting  
\bes
(ks)=\s(s,-1,x_k) \s(k,1,-x_s) \s(s,1,x_k). 
\ees 
Note that the automorphism $(ks)$ of the algebra $A_n$ changes only the positions of the variables $x_k$ and $x_s$.  
Now it is easy to see that  
\bee\label{f3.4}
\s(i,\a,f)^{(ks)}=\s(j,\a,(ks)(f)),
\eee
where $x_j=(ks)(x_i)$. 

\smallskip 

Let $G(A_n)$ be the abstract group with generators (\ref{f3.1}) and defining relations (\ref{f3.2})--(\ref{f3.4}). 

\begin{lm}\label{l3.0}
The subgroup of $G(A_n)$ generated by all 
elements $(ks)$, where $1\leq k\neq s\leq n$, 
is isomorphic to the symmetric group  $S_n$.  
\end{lm}
\Proof
By (\ref{f3.2}) and (\ref{f3.3}), we have 
\bes
(ks)^2=\s(s,-1,x_k) \s(k,1,-x_s) \s(s,1,x_k) \s(s,-1,x_k) \s(k,1,-x_s) \s(s,1,x_k) \\
= \s(s,-1,x_k) \s(k,1,-x_s) \s(s,-1,0) \s(k,1,-x_s) \s(s,1,x_k) \\
= \s(s,-1,x_k) \s(s,-1,0) \s(k,1,-x_s)^{\s(s,-1,0)} \s(k,1,-x_s) \s(s,1,x_k) \\
= \s(s,1,-x_k) \s(k,1,x_s) \s(k,1,-x_s) \s(s,1,x_k) 
= \s(s,1,-x_k) \s(s,1,x_k) =id. 
\ees
Then (\ref{f3.4}) gives 
\bes
(ks)^{(sk)}=\s(s,-1,x_k)^{(sk)} \s(k,1,-x_s)^{(sk)} \s(s,1,x_k)^{(sk)} \\
 =
\s(k,-1,x_s) \s(s,1,-x_k) \s(k,1,x_s)=(sk), 
\ees 
i.e. $(ks)=(sk)$. Now it is not difficult to deduce from (\ref{f3.2})--(\ref{f3.4}) that  
\bes
[(ij),(ks)]=id, \, \, (ik)^{(is)}=(ks), 
\ees
where $i,j,k,s$ are all distinct. It is 
immediate that the given relations imply the defining relations of the group $S_n$ 
with respect to the system of generators $(i\,\,i+1)$, where $1\leq i\leq n-1$, which are indicated in \cite{CM}. 
$\Box$

 By Lemma  \ref{l3.0}, the elements of the symmetric group $S_n$ can be identified with  elements of $G(A_n)$.  
Note that (\ref{f3.4}) can be rewritten as   
\bes
\s(i,\a,f)^{\pi}=\s(\pi^{-1}(i),\a,\pi^{-1}(f)),
\ees
where $\pi \in S_n$.

\smallskip 

It is well known that the group of affine automorphisms $Af_n(F)$ of the algebra $A_n$ is generated by all affine elementary automorphisms.  
\begin{lm}\label{l3.1}
The relations (\ref{f3.2})--(\ref{f3.4}) for elementary affine automorphisms are defining relations of the group $Af_n(F)$. 
\end{lm}
\Proof
Let $\varphi$ be a product of elementary affine automorphisms.  
Suppose that $\varphi=id$. 
 By (\ref{f3.2}) and (\ref{f3.3}), we can represent $\varphi$ in the form  
\bes
\varphi=\s(1,1,\a_1)\s(2,1,\a_2)\ldots \s(n,1,\a_n) \varphi', 
\ees
where $\varphi'$ is a product of elementary linear automorphisms. Obviously, $\a_1=\a_2=\ldots =\a_n=0$. 
Therefore we can assume that $\varphi$ is a product of elementary linear automorphisms.  
 By (\ref{f3.2}) and (\ref{f3.3}), we can easily represent  $\varphi$ in the form  
\bes
\varphi=\s(1,\a_1,0)\s(2,\a_2,0)\ldots \s(n,\a_n,0) \varphi', 
\ees
where $\varphi'$ is a product of elementary automorphisms of the type $\s(i,1,f)$. By (\ref{f3.2})--(\ref{f3.4}), we have   
\bes
\s(k,\a,0)=\s(s,\a,0)^{(ks)} \\
=\s(s,-1,x_k) \s(k,1,-x_s) \s(s,1,x_k) \s(s,\a,0) \s(s,-1,x_k) \s(k,1,-x_s) \s(s,1,x_k) \\
=\s(s,-1,x_k) \s(k,1,-x_s) \s(s,-\a,0) \s(s,1,(1-\a)x_k) \s(k,1,-x_s) \s(s,1,x_k) \\
=\s(s,-1,x_k) \s(s,-\a,0) \s(k,1,\a^{-1}x_s) \s(s,1,(1-\a)x_k) \s(k,1,-x_s) \s(s,1,x_k) \\
=\s(s,\a,0) \s(s,1,-\a x_k) \s(k,1,\a^{-1}x_s) \s(s,1,(1-\a)x_k) \s(k,1,-x_s) \s(s,1,x_k).
\ees
By this relation, we can represent $\varphi$ in the form   
\bes
\varphi=\s(n,\b_n,0) \varphi', 
\ees
where $\varphi'$ is a product of elementary linear automorphisms of the form $\s(i,1,f)$. Hence $\b_n=1$. 
Note that $\s(i,1,f)$ can be represented as a product of automorphisms 
\bee\label{f3.41}
X_{ij}(\lambda)=\s(j,1,\lambda x_i), \ \ \lambda \in F, \ \ i\neq j. 
\eee
Thus, we can assume that $\varphi$ is a product of automorphisms of the form (\ref{f3.41}). 

Let $G$ be the subgroup of $TA_n(F)$ generated by all automorphisms of the form (\ref{f3.41}). 
We define a map 
\bes
J : G \longrightarrow SL_n(F), 
\ees
where $J(\psi)$ is the Jacobian matrix of $\psi\in G$. It is easy to check that 
\bes
J(X_{ij}(\lambda))=E_{ij}(\lambda), 
\ees 
and that $J$ is an isomorphism of the groups. 

Now it is sufficient to prove that every relation of the group $SL_n(F)$ is a corollary of (\ref{f3.2})--(\ref{f3.4}). 
Obviously,  (\ref{f3.2})--(\ref{f3.3}) cover the Steinberg relations (see, for example \cite{Milnor}). 
Besides, according to \cite{Milnor}, we need to check the relation  
\bes
\{u,v\}=id, \ \ 0\neq u,v\in F, 
\ees 
where 
\bes
\{u,v\}=h_{ij}(uv) h_{ij}(u)^{-1} h_{ij}(v)^{-1}, \\
h_{ij}(u)=w_{ij}(u) w_{ij}(-1), \\
w_{ij}(u)=X_{ij}(u) X_{ji}(-u^{-1}) X_{ij}(u). 
\ees
Applying (\ref{f3.2})--(\ref{f3.4}) we have 
\bes
w_{ij}(u)= \s(j,1,u x_i) \s(i,1,-u^{-1} x_j) \s(j,1,u x_i) \\
= \s(j,1,u x_i) \s(i,u,0) \s(i,1,-x_j) \s(i,u^{-1},0) \s(j,1,u x_i) \\
= \s(i,u,0) \s(j,1,u x_i)^{\s(i,u,0)} \s(i,1,-x_j) \s(j,1,u x_i)^{\s(i,u,0)} \s(i,u^{-1},0) \\
= \s(i,u,0) \s(j,1,x_i) \s(i,1,-x_j) \s(j,1,x_i) \s(i,u^{-1},0) \\
= \s(i,u,0) \s(j,-1,0) \s(j,-1,x_i) \s(i,1,-x_j) \s(j,1,x_i) \s(i,u^{-1},0) \\
= \s(i,u,0) \s(j,-1,0) (ij) \s(i,u^{-1},0) = (ij) \s(i,u,0)^{(ij)} \s(j,-1,0)^{(ij)} \s(i,u^{-1},0)\\
= (ij) \s(j,u,0) \s(i,-1,0) \s(i,u^{-1},0)= (ij) \s(j,u,0) \s(i,-u^{-1},0). 
\ees
Consequently, 
\bes
h_{ij}(u)=w_{ij}(u) w_{ij}(-1) = (ij) \s(j,u,0) \s(i,-u^{-1},0) (ij) \s(j,-1,0) \s(i,1,0) \\
= \s(j,u,0)^{(ij)} \s(i,-u^{-1},0)^{(ij)} \s(j,-1,0)  \\ 
= \s(i,u,0) \s(j,-u^{-1},0) \s(j,-1,0) = \s(i,u,0) \s(j,u^{-1},0).  
\ees
Hence  
\bes
\{u,v\}=h_{ij}(uv) h_{ij}(u)^{-1} h_{ij}(v)^{-1} \\
=  \s(i,uv,0) \s(j,(uv)^{-1},0) \s(i,u,0) \s(j,u^{-1},0) \s(i,v,0) \s(j,v^{-1},0)=id.
\ees
 Thus we can say that every relation of the group $SL_n(F)$ follows from (\ref{f3.2})--(\ref{f3.4}). 
$\Box$

Assume that   
\bee\label{f3.5}
\te=\phi_1\phi_2\ldots \phi_r \in TA_n(F), 
\eee
where $\phi_i$, $1\leq i\leq r$, are elementary automorphisms. 
Put  
\bes
\psi_i=\phi_1\phi_2\ldots\phi_i, \,\, 0\leq i\leq r.  
\ees
In particular, we have 
\bes
\psi_r=\te,\,\,\psi_0=id.  
\ees
To (\ref{f3.5}) corresponds the sequence of elementary transformations  
\bee\label{f3.6}
id=\psi_0\stackrel{\phi_1}{\rightarrow}\psi_1\stackrel{\phi_2}{\rightarrow}\psi_2\stackrel{\phi_3}{\rightarrow}\ldots \stackrel{\phi_r}{\rightarrow}\psi_r=\te.  
\eee
Note that the representations (\ref{f3.5}) and (\ref{f3.6}) of the automorphism  $\te$ are equivalent. 
If (\ref{f3.6}) is a minimal representation of $\te$, then the representation (\ref{f3.5}) will be also called a {\em minimal representation} of $\te$.
\begin{theor}\label{t3.1} Let $F$ be a field of characteristic $0$. 
Then the relations (\ref{f3.2})--(\ref{f3.4}) are defining relations of the group $TA_3(F)$
with respect to the generators (\ref{f3.1}). 
\end{theor}
{\bf Plan of the proof.}
Assume that   
\bee\label{f3.7}
\varphi_1\varphi_2\ldots\varphi_k=id=(x_1,x_2,x_3),
\eee
where $\varphi_i$, $1\leq i\leq k$, are elementary automorphisms. Put    
\bes
\te_i=\varphi_1\varphi_2\ldots\varphi_i,\,\,0\leq i\leq k.   
\ees
In particular, we have $\te_0=\te_k=(x_1,x_2,x_3)$. To (\ref{f3.7}) corresponds the sequence of elementary transformations  
\bee\label{f3.8}
id=\te_0\stackrel{\varphi_1}{\rightarrow} \te_1\stackrel{\varphi_2}{\rightarrow} \ldots \stackrel{\varphi_k}{\rightarrow} \te_k=id.
\eee

Denote by $d=max\{\deg\,\te_i | 0\leq i\leq k\}$ the {\em width}  
of the sequence (\ref{f3.8}). Let $i_1$ be the minimal number and $i_2$ be the maximal 
number which satisfy the equations $\deg\,\te_{i_1}=d$ and $\deg\,\te_{i_2}=d$. Put $q=i_2-i_1$. 
The pair $(d,q)$ will be called the {\em exponent} of the relation (\ref{f3.7}). 

To prove the theorem, we show that (\ref{f3.7}) 
 follows from (\ref{f3.2})--(\ref{f3.4}). Assume that our theorem is not true. 
 Call a relation of the form (\ref{f3.7})  
  {\em trivial} if it follows from (\ref{f3.2})--(\ref{f3.4}). 
 We choose a nontrivial relation (\ref{f3.7})  with the minimal exponent $(d,q)$ 
 with respect to the lexicographic order. To arrive at a contradiction, 
 we show that (\ref{f3.7}) is also trivial.  

If $d=3$, then Lemma \ref{l3.1} gives the triviality of the relation (\ref{f3.7}). 
Therefore we can assume that $d>3$.

Our plan is to change the product  (\ref{f3.7}) by using (\ref{f3.2})--(\ref{f3.4}) and 
to obtain a new sequence (\ref{f3.8}) whose exponent is strictly less than $(d,q)$. 
The proof of the theorem will be completed by Lemmas \ref{l3.2}--\ref{l3.9}. $\Box$

\medskip

 Denote by $t=[\frac{q}{2}]$ the integral part of $\frac{q}{2}$.
Put also   
\bes
\phi=\te_{i_1+t-1},\,\, 
\te=\te_{i_1+t},\,\, \tau=\te_{i_1+t+1}.
\ees
 Then we have  
\bee\label{f3.9}
\phi \stackrel{\varphi_{i_1+t}}{\longrightarrow} \theta \stackrel{\varphi_{i_1+t+1}}{\longrightarrow} \tau. 
\eee

\begin{lm}\label{l3.2} 
The following statements are true:  
\begin{itemize}
\item[(1)]  $d$ is the width of $\te$, $t$ is the height of $\te$, and  
\bee\label{f3.10}
\te=\varphi_1\varphi_2\ldots \varphi_{i_1+t} 
\eee  
is a minimal representation of $\te$. 
\item[(2)]  If $q$ is an even number, then  
\bes
\te=\varphi_{k}^{-1}\varphi_{k-1}^{-1}\ldots \varphi_{i_1+t+1}^{-1}  
\ees
is also a minimal representation of $\te$. 
\item[(3)]  If $q$ is an odd number, then $(d(\tau),t(\tau))=(d,t)$ and  
\bes
\tau=\varphi_{k}^{-1}\varphi_{k-1}^{-1}\ldots \varphi_{i_1+t+2}^{-1}  
\ees
is a minimal representation of $\tau$. Moreover, in (\ref{f3.7}) 
the product (\ref{f3.10}) can be replaced by an arbitrary minimal representation of 
 $\te$. 
\end{itemize} 
\end{lm}
\Proof
Assume that $(d(\te),t(\te))<(d,t)$ and let (\ref{f3.5}) be a minimal representation of $\te$. Then (\ref{f3.7}) is a consequence of the equalities    
\bee\label{f3.11}
\varphi_1\varphi_2\ldots \varphi_{i_1+t}\phi_r^{-1}\ldots \phi_2^{-1}\phi_1^{-1} =id,  
\eee  
\bee\label{f3.12}
\phi_1\phi_2\ldots\phi_r\varphi_{i_1+t+1}\ldots \varphi_{k-1}\varphi_k=id.  
\eee  
To (\ref{f3.11}) corresponds the sequence of elementary transformations  
\bes
(x_1,x_2,x_3)\rightarrow \te_1\rightarrow\ldots\rightarrow\te_{i_1+t}=
\te=\psi_r\rightarrow\psi_{r-1}\ldots\rightarrow\psi_1\rightarrow(x_1,x_2,x_3), 
\ees
and to (\ref{f3.12}) corresponds  
\bes
(x_1,x_2,x_3)\rightarrow\psi_1\rightarrow\ldots\rightarrow\psi_r=\te=\te_{i_1+t}\rightarrow\te_{i_1+t+1}
\rightarrow\ldots\rightarrow\te_{k-1}\rightarrow(x_1,x_2,x_3). 
\ees
Since $(d(\te),t(\te))<(d,t)$, it follows that (\ref{f3.11}) and (\ref{f3.12}) 
have exponents strictly less than $(d,q)$. This gives the first statement of the lemma. 

It is obvious that  
\bes
\varphi_k^{-1}\varphi_{k-1}^{-1}\ldots\varphi_1^{-1}=id 
\ees
has the same exponent $(d,q)$. Applying the first statement of the lemma  to this relation, 
we get the second statement of the lemma, as well as the minimality of the representation of $\tau$  
if $q$ is an odd number. If $q$ is an odd number, then (\ref{f3.11}) 
has exponent strictly less than $(d,q)$, and  (\ref{f3.12}) has the exponent $(d,q)$. 
Consequently, (\ref{f3.7}) and (\ref{f3.12}) are equivalent modulo  
(\ref{f3.2})--(\ref{f3.4}). Thus $\te$ can be changed by 
an arbitrary minimal representation in (\ref{f3.12}). $\Box$

\begin{lm}\label{l3.3} 
If $t\geq 1$, then the relation (\ref{f3.7}) is trivial.  
\end{lm}
\Proof
Our aim is to change the product 
$\varphi_{i_1+t}\varphi_{i_1+t+1}$ by (\ref{f3.2})--(\ref{f3.4}) and to get a new sequence 
(\ref{f3.8}) whose exponent is strictly less than $(d,q)$, i.e. to show that (\ref{f3.7}) is trivial.  

Since $t\geq 1$, according to Proposition \ref{p2.1} $\te$ admits a reduction of  
types  I--IV. We restrict ourselves only to the case when $\te$ admits a reduction of 
 type IV. The other cases can be considered similarly.   

Let $\te=(f_1,f_2,f_3)$ satisfy the conditions of Definition \ref{d1.4}. 
According to Lemma \ref{l3.2}, the representation (\ref{f3.10}) of $\te$ is minimal. 
Then $\phi$ can be calculated by using Proposition \ref{p2.1}. 

{\bf Case I: $q$ is even, $q=2t$.} 
According to Lemma \ref{l3.2} $\tau$ can also be calculated by using Proposition \ref{p2.1}. 

Assume that $t=2$. By Proposition \ref{p2.1} the automorphisms $\phi$ and $\tau$ have the same form, i.e.  either 
\bes
\phi=(f_1,\eta_1g_2+\kappa_1g_1+\nu_1,f_3),\,\,\, \tau=(f_1,\eta_2g_2+\kappa_2g_1+\nu_2,f_3),
\ees
 or 
\bes
\phi=(\rho_1g_1+\s_1,f_2,f_3),\,\,\, \tau=(\rho_2g_1+\s_2,f_2,f_3). 
\ees
 By (\ref{f3.2}), in both cases $\varphi_{i_1+t}\varphi_{i_1+t+1}$ can be replaced by an elementary automorphism. 
Obviously, (\ref{f3.7}) will then be replaced  by a relation whose exponent is strictly less than $(d,q)$. 

Assume that $t=3$. According to Proposition \ref{p2.1}, the automorphisms  
 $\phi$ and  $\tau$ have the forms (\ref{f2.3}) and (\ref{f2.4}). If $\phi$ and $\tau$
 have the same form (\ref{f2.3}) (or (\ref{f2.4})), then as above,  by (\ref{f3.2}), we obtain the triviality of (\ref{f3.7}). 
 Suppose that $\phi$ has the form (\ref{f2.3}), and  $\tau$ has the form (\ref{f2.4}). It is immediate that  
\bes
\varphi_{i_1+t}=\s(2,\eta^{-1},F(x_1,x_3)), \ \ \varphi_{i_1+t+1}=\s(1,\rho,G(x_3)). 
\ees
By (\ref{f3.3}) we get 
\bes
\varphi_{i_1+t}\varphi_{i_1+t+1}=\s(1,\rho,G(x_3))\s(2,\eta^{-1},F_1(x_1,x_3)).
\ees
 We replace $\varphi_{i_1+t}\varphi_{i_1+t+1}$ in (\ref{f3.7}) according to this equality.   
Then $\te$ in (\ref{f3.8}) can be changed to     
\bes
\te'=(\rho g_1+\a,\eta g_2+\kappa g_1+\nu,f_3).
\ees
Note that $\deg\,\te'=\frac{13n}{2}$ and after such replacement the exponent 
$(d,q)$ of the sequence (\ref{f3.8}) does not change. As before $\te'$ admits a 
reduction of type IV. But according to Proposition \ref{p2.1}, in this case  
we have $t(\te')=1$. This contradicts Lemma \ref{l3.2}, since $t(\te')<t$. 

If $t=1$, then according to Proposition \ref{p2.1} we have  
\bes
\phi=(g_1,g_2,\rho_1g_3+\a_1),\,\,\,\tau=(g_1,g_2,\rho_2g_3+\a_2).
\ees
By (\ref{f3.2}), we obtain the triviality of (\ref{f3.7}) as above. 

{\bf Case II: $q$ is odd, $q=2t+1$.} 
If $(\a_1,\a_2)=0$ and $(\b_1,\b_2,\b_3,\b_4)=0$ in the conditions of Definition \ref{d1.4}, 
then $\te=(g_1,g_2,f_3)$. Moreover, according to Proposition \ref{p2.1}, we have $t=1$ and 
\bes
\phi=(g_1,g_2,\rho g_3+\a). 
\ees
Assume that  
\bee\label{f3.13}
\tau=(\lambda f_1+a,f_2,f_3),\,\,\,a\in<f_2,f_3>.
\eee
By Lemma \ref{l3.2} we have $(d(\tau),t(\tau))=(d,t)$. Since  
 $d=d(\te)=\frac{13n}{2}$, from here we get $\deg\,\tau \leq \frac{13n}{2}$. 
 With regard to the inequalities of Proposition \ref{p1.1}, we have $\deg\,(\lambda f_1+a)<\frac{5n}{2}$ and $\deg\,a<\frac{5n}{2}$. 
 Applying Lemma \ref{l1.4} gives also  
\bes
a=\a_0'+\a_1'f_3+\a_2'f_3^2.
\ees
Consequently,  
\bes
\lambda f_1+a=\lambda g_1+\a_0'+\a_1'f_3+\a_2'f_3^2.
\ees
If $(\a_1',\a_2')\neq 0$, then $\tau$ admits a reduction of  type IV. Moreover, $t(\tau)>1$, which contradicts the equality  
$(d(\tau),t(\tau))=(d,t)$. Consequently,  
\bes
\tau=(\lambda g_1+\a_0',g_2,f_3).
\ees
We have  
\bes
\varphi_{i_1+t}=\s(3,\rho^{-1},F(x_1,x_2)),\,\,\,\varphi_{i_1+t+1}=\s(1,\lambda,\a_0').
\ees
According to (\ref{f3.3}), we obtain  
\bes
\varphi_{i_1+t}\varphi_{i_1+t+1}=\s(1,\lambda,\a_0')\s(3,\rho^{-1},F_1(x_1,x_2)).
\ees
After the corresponding replacement, instead of $\te$ we have   
\bes
\te'=(\lambda g_1+\a_0',g_2,\rho g_3+\a). 
\ees
Note that $\te'$ admits an essential elementary reduction, i.e. 
$t(\te')=0$. 

Now assume that  
\bee\label{f3.14}
\tau=(f_1,\eta f_2+a,f_3),\,\,\,a\in <f_1,f_3>. 
\eee
Using the same arguments as above we get  
\bes
\tau=(g_1,\eta g_2+\kappa g_1+\nu,f_3). 
\ees
We have  
\bes
\varphi_{i_1+t}=\s(3,\rho^{-1},F(x_1,x_2)),\,\,\,\varphi_{i_1+t+1}=\s(2,\eta,\kappa x_1+\nu).
\ees
The relation (\ref{f3.3}) gives 
\bes
\varphi_{i_1+t}\varphi_{i_1+t+1}=\s(2,\eta,\kappa x_1+\nu)\s(3,\rho^{-1},F_1(x_1,x_2)).
\ees
After such replacement, instead of $\te$ we have  
\bes
\te'=(g_1,\eta g_2+\kappa g_1+\nu,\rho g_3+\a),  
\ees
and this automorphism also admits an essential elementary reduction.  

If the elementary reduction $\te\rightarrow \tau$ replaces the element 
 $f_3$ of the automorphism $\te$, then applying (\ref{f3.2}) also gives the triviality of (\ref{f3.7}). 

We now consider the case when $(\a_1,\a_2)=0$ and $(\b_1,\b_2,\b_3,\b_4)\neq 0$. 
Then according to Proposition \ref{p2.1} we have $t=2$ and  $\phi$ has the form (\ref{f2.3}), i.e.  
\bes
\phi=(g_1,\eta g_2+\kappa g_1+\nu,f_3).   
\ees
If $\tau$ has the form (\ref{f3.14}), then (\ref{f3.2}) gives the triviality of 
 (\ref{f3.7}). Assume that $\tau$ has the form (\ref{f3.13}). Then by the same discussion as above we get  
\bes
\tau=(\rho g_1+\a,f_2,f_3).
\ees
We have  
\bes
\varphi_{i_1+t}=\s(2,\eta^{-1},F(x_1,x_2)),\,\,\,\varphi_{i_1+t+1}=\s(1,\rho,\a).
\ees
By (\ref{f3.3}), we get  
\bes
\varphi_{i_1+t}\varphi_{i_1+t+1}=\s(1,\rho,\a)\s(3,\eta^{-1},F_1(x_1,x_2)).
\ees
After the corresponding replacement, instead of $\te$ we obtain  
\bes
\te'=(\rho g_1+\a,\eta g_2+\kappa g_1+\nu,f_3).
\ees
Proposition \ref{p2.1} gives $t(\te')=1<t$; a contradiction. 

Assume that  
\bee\label{f3.15}
\tau=(f_1,f_2,\rho f_3+a),\,\,\,a\in <f_1,f_2>.
\eee
By Proposition \ref{p1.1} we have $\deg\,f_1+\deg\,f_2>\frac{9n}{2}$. 
Since $d(\te)=\frac{13n}{2}\geq \deg\,\tau$, we have $\deg\,(\rho f_3+a)<2n$ and $\deg\,a<2n$. 
Lemma \ref{l1.4} gives $a=\a \in F$. Then $\varphi_{i_1+t+1}=\s(3,\rho,\a)$. Hence   
\bes
\varphi_{i_1+t}\varphi_{i_1+t+1}=\s(3,\rho,\a)\s(2,\eta^{-1},F_1(x_1,x_2)).
\ees
After this replacement, $\te$ is changed to   
\bes
\te'=(g_1,\eta g_2+\kappa g_1+\nu,\rho f_3+a). 
\ees
Proposition \ref{p2.1} gives $t(\te')=1$, which also leads to a contradiction.  

The case when $(\a_1,\a_2)\neq 0$ and $(\b_1,\b_2,\b_3,\b_4)=0$ can be considered analogously.  

Assume that $(\a_1,\a_2)\neq 0$ and $(\b_1,\b_2,\b_3,\b_4)\neq 0$. 
According to Proposition \ref{p2.1}, we have $t=3$. Now we use the statement (3) of Lemma \ref{l3.2} about the 
arbitrariness of the minimal representation of 
 $\te$ in (\ref{f3.7}). If $\tau$ has the form (\ref{f3.13}), 
then we can assume that $\phi$ has the form (\ref{f2.4}), and if $\tau$ has the form (\ref{f3.14}), 
then we can assume that $\phi$ has the form (\ref{f2.3}). 
By (\ref{f3.2}), in both cases we can decrease the exponent of (\ref{f3.7}). 
If $\tau$ has the form (\ref{f3.15}), then using the same arguments we get  
\bes
\tau=(f_1,f_2,\rho f_3+\a).
\ees
Assume that $\phi=(\lambda g_1+\b,f_2,f_3)$. Applying  
 (\ref{f3.3}) to $\varphi_{i_1+t}\varphi_{i_1+t+1}$ changes $\te$ to     
\bes
\te'=(g_1,g_2,\rho f_3+a). 
\ees
Then we have $t(\te')=2$, which also leads to a contradiction. $\Box$

Now we begin to consider the most complicated case when $t=0$, i.e. $q=0, 1$. Put $\te=(f_1,f_2,f_3)$. According to Lemma \ref{l3.2},    
 $t=0$ is the height of $\te$, i.e. $\te$ admits an essential elementary reduction. Moreover, $\phi$ is an essential elementary reduction of $\te$. 
 If $q=0$, then $\tau$ is also an esssential elementary reduction of $\te$. 
If $q=1$, then we have $(d(\tau),t(\tau))=(d,t)$. Consequently,  $\deg\,\tau\leq \deg\,\te=d$. Thus we can assume that  
\bee\label{f3.16}
\t=(f_1,f_2,f),\ \  f=\a f_3+a,\ \  a=a(f_1,f_2)\in \langle f_1,f_2 \rangle,\ \
\deg a\leq \deg f_3 .
\eee

\begin{lm}\label{l3.4}
If $\phi$ reduces the element $f_3$ of $\te$, then the relation (\ref{f3.7}) is trivial.  
\end{lm}
\Proof
Applying (\ref{f3.2}) we can replace $\varphi_{i_1+t}\varphi_{i_1+t+1}$ 
 by an elementary automorphism.  Obviously, this replacement also decreases the exponent of (\ref{f3.8}).
$\Box$

By Lemma \ref{l3.4}, we can assume that $\phi$ reduces one of the elements  $f_1$ and $f_2$ of $\te$. 
 Without loss of generality, later on we consider that  $\phi$ reduces the element  $f_2$ of $\te$, 
i.e. $\phi=(f_1,g_2,f_3)$ and $\deg\,g_2<\deg\,f_2$.   
\begin{lm}\label{l3.5}
If $\phi'$ reduces the element $f_2$ of $\te$, then in 
 (\ref{f3.9}) the automorphism $\phi$ can be replaced by $\phi'$. 
\end{lm}
\Proof
According to (\ref{f3.2}), in this case the elementary 
transformation $\phi\rightarrow\te$ can be changed to 
 $\phi\rightarrow\phi'\rightarrow\te$. Since $\deg\,\phi'<\deg\,\te=d$, 
 the exponent $(d,q)$ of the sequence (\ref{f3.8}) does not change after this replacement. 
 But in the new sequence  (\ref{f3.8}) we have $\phi'$  instead of $\phi$. $\Box$

Taking  this lemma into account, we can assume that  
\bee\label{f3.17}
\phi=(f_1,g_2,f_3),\ \ g_2=f_2-b,\ \  b\in \langle f_1,f_3 \rangle,\ \  \deg g_2<\deg f_2.
\eee

\begin{lm}\label{l3.6}
If one of the following conditions is fulfilled, then (\ref{f3.7}) is trivial: 

$(1)$ $\overline{f_2}\in \langle\overline{f_1}\rangle$;

$(2)$ $\overline{f_3}\in \langle\overline{f_1}\rangle$;

$(3)$ $a$ does not depend on $f_2$;

$(4)$ $\overline{f_2}=\b \overline{f_3}+\gamma \overline{f_1^k}$; 

$(5)$ $\overline{f_1}$ and $\overline{f_3}$ are algebraically independent. 
\end{lm}

\Proof
Assume that $\overline{f_2}\in \langle\overline{f_1}\rangle$ and $\overline{f_2}=\b \overline{f_1}^l$. According to Lemma \ref{l3.5}, we can suppose that 
 $g_2=f_2-\b f_1^l$. Then  
\bes
\varphi_{i_1+t}=\s(2,1,\b x_1^l),\ \ \varphi_{i_1+t+1}=\s(3,\a,a(x_1,x_2)). 
\ees
By (\ref{f3.3}), we have 
\bes
\varphi_{i_1+t}\varphi_{i_1+t+1}=\s(3,\a,a_1(x_1,x_2))\s(2,1,\b x_1^l). 
\ees
After the corresponding replacement in (\ref{f3.7}), 
 $\te$ is replaced by $\te'=(f_1,g_2,f)$ in (\ref{f3.8}). 
Since  $\deg\,f_1+\deg\,g_2+\deg\,f<d$, the exponent of (\ref{f3.7}) is decreased.  

\smallskip

Assume that $\overline{f_3}\in \langle\overline{f_1}\rangle$ and $\overline{f_3}=T(\overline{f_1})$. Put $g_3=f_3-T(f_1)$. 
According to (\ref{f3.3}), we have  
\bes
\varphi_{i_1+t}=\s(2,1,b(x_1,x_3))=\s(3,1,-T(x_1))\s(2,1,b_1(x_1,x_3))\s(3,1,T(x_1)). 
\ees
After the corresponding replacement in (\ref{f3.7}), the elementary transformation  
$\phi\rightarrow\te$ is replaced by the sequence of elementary transformations  
\bes
\phi=(f_1,g_2,f_3)\rightarrow (f_1,g_2,g_3)\rightarrow (f_1,f_2,g_3)
\rightarrow (f_1,f_2,f_3)=\te.
\ees
Since $d(\phi), \deg\,(f_1,g_2,g_3), \deg\,(f_1,f_2,g_3)<d=\deg\,\te$, 
 the new sequence (\ref{f3.8}) has the same exponent $(d,q)$. 
However, instead of $\phi$ we have $(f_1,f_2,g_3)$. 
By Lemma \ref{l3.4}, we obtain the triviality of (\ref{f3.7}). 

\smallskip

Assume that $a$ does not depend on $f_2$. By (\ref{f3.3}) we have  
\bes
\varphi_{i_1+t}\varphi_{i_1+t+1}=\s(2,1,b(x_1,x_3))\s(3,\a,a(x_1))=\s(3,\a,a(x_1))\s(2,1,b_1(x_1,x_3)). 
\ees
After the corresponding replacement in (\ref{f3.7}), instead of $\te$ we obtain $\te'=(f_1,g_2,f)$. 
Since $\deg\,(f_1,g_2,f)<d$, this replacement also decreases the exponent of (\ref{f3.7}). 

\smallskip

Consider the case when $\overline{f_2}=\b \overline{f_3}+\gamma \overline{f_1^k}$. 
By Lemma \ref{l3.6}(1) we can assume that $\b\neq 0$. By Lemma \ref{l3.5} we can also 
assume that $b=\b f_3+\gamma f_1^k$. Consequently,  
\bes
g_2&=&f_2-\b f_3-\gamma f_1^k,\\
f_2&=&g_2+\b f_3+\gamma f_1^k,\\
f_3&=&-\frac{1}{\b} g_2+\frac{1}{\b}f_2-\frac{\gamma}{\b} f_1^k. 
\ees
These equalities justify the sequence of elementary transformations 
\bes
(f_1,f_3,g_2)\rightarrow (f_1,f_2,g_2)\rightarrow (f_1,f_2,f_3)=\te.
\ees
We have  
\bes
\varphi_{i_1+t}=\s(2,1,\b x_3+\gamma x_1^k)). 
\ees
Applying  (\ref{f3.2}) and (\ref{f3.3}) we get  
\bes
\varphi_{i_1+t}=\s(2,1,\b x_3)\s(2,1,\gamma x_1^k)\\
=\s(2,1,\b x_3)\s(3,-\b,x_2)\s(3,-\frac{1}{\b},\frac{1}{\b}x_2)\s(2,1,\gamma x_1^k)\\
=\s(2,1,\b x_3)\s(3,-\b,x_2)\s(2,1,\gamma x_1^k)\s(3,-\frac{1}{\b},\frac{1}{\b}x_2)^{\s(2,1,\gamma x_1^k)}\\
=\s(2,1,\b x_3)\s(3,-\b,x_2)\s(2,1,\gamma x_1^k)\s(3,-\frac{1}{\b},\frac{1}{\b}(x_2-\gamma x_1^k))\\
=\s(2,1,\b x_3)\s(3,-\b,x_2)\s(2,\frac{1}{\b},-\frac{1}{\b}x_3)\s(2,\b,x_3+\gamma x_1^k)\s(3,-\frac{1}{\b},\frac{1}{\b}(x_2-\gamma x_1^k)).
\ees
Since the transposition $(23)\in S_3$ has a linear representation  
\bes
(23)=\s(2,1,\b x_3)\s(3,-\b,x_2)\s(2,\frac{1}{\b},-\frac{1}{\b}x_3), 
\ees
we obtain   
\bes
\varphi_{i_1+t}=(23)\s(2,\b,x_3+\gamma x_1^k))\s(3,-\frac{1}{\b},\frac{1}{\b}(x_2-\gamma x_1^k)).
\ees
Then  
\bee\label{f3.18}
\te=(23)\varphi_1^{(23)}\varphi_2^{(23)}\ldots\varphi_{i_1+t-1}^{(23)}\s(2,\b,x_3+\gamma x_1^k))\s(3,-\frac{1}{\b},\frac{1}{\b}(x_2-\gamma x_1^k)),
\eee
where $\varphi_i^{(23)}$ are elementary automorphisms, according to (\ref{f3.4}). To (\ref{f3.18}) corresponds the sequence of elementary transformations  
\bes
(x_1,x_2,x_3)\mapsto(x_1,x_3,x_2)\rightarrow\te_1'\rightarrow\te_2'\rightarrow\ldots\rightarrow\te_{i_1+t-1}'\\ 
=(f_1,f_3,g_2)\rightarrow(f_1,f_2,g_2)\rightarrow(f_1,f_2,f_3)=\te, 
\ees
where $\te_i'$ is obtained from $\te_i$ only by the permutation of the second and the 
third components, and the automorphism $(x_1,x_2,x_3)\mapsto(x_1,x_3,x_2)$ is a composition of three elementary 
linear transformations. 

If in (\ref{f3.7}) we replace $\te$ by (\ref{f3.18}), 
then the exponent of (\ref{f3.8}) remains the same. But instead of  
$\phi$ we have $(f_1,f_2,g_2)$, and Lemma \ref{l3.4} gives the triviality of (\ref{f3.7}). 

\smallskip

We now consider the case when $\overline{f_1}$ and $\overline{f_3}$ are algebraically independent. 
Then $\overline{f_2}=\overline{b}\in \langle\overline{f_1},\overline{f_3}\rangle$. By Lemma \ref{l3.6}(1) we can assume that 
$\overline{f_2}\notin \langle\overline{f_1}\rangle$, 
i.e. $\overline{f_2}$ depends on $\overline{f_3}$. Consequently, 
$\deg f_3 \leq \deg f_2$. If $\overline{f_1}$ and $\overline{f_2}$ are algebraically dependent, then it follows that  
$\overline{f_1}$ and $\overline{f_3}$ are algebraically dependent.  
Consequently, $\overline{f_1}$ and $\overline{f_2}$ are algebraically independent. 
Then $\overline{a}\in \langle \overline{f_1},\overline{f_2}\rangle$. By Lemma \ref{l3.6}(3) we can assume that $a$ necessarily contains $f_2$. 
Then $\deg\,f_2\leq \deg\,a\leq \deg\,f_3 $, i.e.  $\deg\,f_2=\deg\,f_3$. 
 Hence  
\bes
\overline{b}=\overline{f_2}=\b \overline{f_3}+\gamma \overline{f_1^k}.
\ees
From the statement (4) of the lemma we obtain that (\ref{f3.7}) is trivial. 
$\Box$

\smallskip

So, by Lemma \ref{l3.6}, we can assume that $\overline{f_1}$ and $\overline{f_3}$ are algebraically 
dependent and that $\overline{f_3}\notin \langle\overline{f_1}\rangle$. 
It remains to consider the following three cases separately: 

(1) $f_1,f_3$ is a $\ast$-reduced pair and $\deg f_1<\deg f_3$;

(2) $f_1,f_3$ is a $\ast$-reduced pair and $\deg f_3 <\deg f_1$;

(3) $\overline{f_1}\in \langle\overline{f_3}\rangle$ and $\deg f_1>\deg f_3$.

\begin{lm}\label{l3.7}
If $f_1, f_3$ is a $\ast$-reduced pair and $\deg f_1<\deg f_3 $,
then (\ref{f3.7}) is trivial. 
\end{lm}
\Proof 
We first consider the case when $\deg f_2>\deg f_3 $. If $\overline{f_1}$ and $\overline{f_2}$ are algebraically independent, then 
$\overline{a}\in<\overline{f_1},\overline{f_2}>$. Since    
$\deg\,a\leq\deg\,f_3<\deg\,f_2$, we have $a\in \langle f_1\rangle$. 
 Lemma \ref{l3.6} gives the triviality of (\ref{f3.7}). 

Suppose that $\overline{f_1}$ and $\overline{f_2}$ are algebraically dependent. 
By Lemma \ref{l3.6}, we can assume that $\overline{f_2}\notin \langle\overline{f_1}\rangle$. Since  
 $\deg\,f_1<\deg\,f_3<\deg\,f_2$, therefore $f_1,f_2$ is a $\ast$-reduced pair.  

Assume that $\deg\,a<N=N(f_1,f_2)$. Since $\deg\,a\leq \deg\,f_3<\deg\,f_2$, 
 Corollary \ref{c1.1} gives $a\in<f_1>$. By Lemma \ref{l3.6}, the relation (\ref{f3.7}) is trivial.

Therefore we can assume that $\deg\,a\geq N$. Then $\deg\,f_2>N=N(f_1,f_2)$.  
By the definition of $N=N(f_1,f_2)$, this is possible only if $f_1,f_2$ 
is a 2-reduced pair. Consequently,   $\deg\,f_1=2n$, $\deg\,f_2=sn$, where $s\geq 3$ is an odd number, and  moreover  
\bes
N=n(s-2)+\deg\,[f_1,f_2]\leq \deg\,a\leq \deg\,f_3.
\ees 

Let $(f_1,f_2,g_3)$ be an elementary reduction of $\te$ such that $g_3$ is an irreducible 
element of $(f_1,f_2,g_3)$. If $f_3$ is an irreducible element of $\te$, then we put $g_3=f_3$. 
Consequently, $(f_1,f_2,g_3)$ satisfies the conditions of Proposition \ref{p1.2}. 
Assume that Proposition \ref{p1.2}(1) is valid for $(f_1,f_2,g_3)$, i.e.    
\bes
\deg\,g_3<n(s-2)+\deg\,[f_1,f_2]=N.
\ees 
Since $\deg\,f_3\geq N$, it means that  
\bes
f_3=g_3+c(f_1,f_2), \ \ \overline{c(f_1,f_2)}=\overline{f_3}. 
\ees
By Lemma \ref{l3.6} we can assume that $\overline{c(f_1,f_2)}=\overline{f_3}\notin<f_1>$. Since $\deg\,c(f_1,f_2)=\deg\,f_3<\deg\,f_2$, 
we have $\overline{c(f_1,f_2)}\notin<\overline{f_1},\overline{f_2}>$. 
It is easy to deduce from (\ref{f1.1}) that $\deg_yc(x,y)$ equals $2$ or $4$. 
Then $\deg_y\frac{\partial c}{\partial y}(x,y)$ equals $1$ or $3$, and (\ref{f1.1}) gives $\frac{\partial c}{\partial y}(f_1,f_2)\geq \deg\,f_2$. We have   
\bes
\deg\,[f_1,g_3]\leq \deg\,f_1+\deg\,g_3<ns+\deg\,[f_1,f_2]=\deg\,f_2+\deg\,[f_1,f_2]. 
\ees 
Then   
\bes
[f_1,f_3]=[f_1,g_3]+[f_1,f_2]\frac{\partial c}{\partial y}(f_1,f_2)
\ees
gives that  
\bes
\deg\,[f_1,f_3]=\deg\,[f_1,f_2]\frac{\partial c}{\partial y}(f_1,f_2)\geq\deg\,[f_1,f_2]+\deg\,f_2. 
\ees
Consequently, $N(f_1,f_3)>\deg\,f_2$. By (\ref{f3.17}) we have  $\deg\,b=\deg\,f_2<N(f_1,f_3)$. 
Corollary \ref{c1.1} gives $\overline{f_2}\in<\overline{f_1},\overline{f_3}>$. Since $\deg\,f_1f_3, \deg\,f_3^2 >\deg\,f_2$, we obtain    
$\overline{f_2}\in<\overline{f_1}>$; a contradiction.  

We now assume that one of the statements (2) and  (3) of Proposition \ref{p1.2} is valid for $(f_1,f_2,g_3)$. Combining both the cases we can say that  
\bes
f_3=h_3+r(f_1,f_2),\,\, \deg\,f_2=3n,\,\, \deg\,(f_2-\b h_3^2)\leq 2n,\,\, \deg\,[f_1,h_3]<n+\deg\,[f_1,f_2]. 
\ees
Since $\deg\,f_3<2n$ and $\deg\,h_3=\frac{3n}{2}$, it follows that $\deg\,r(f_1,f_2)<2n$. 
By Corollary \ref{c1.2}, we have $r(f_1,f_2)=\lambda w(f_1,f_2)+\mu$, where $w(x,y)$ is a derivative polynomial of the pair $f_1,f_2$. 
 If $\lambda\neq 0$, then an immediate calculation gives $\deg\,[f_1,f_3]>\deg\,f_2$.  
 This leads to a contradiction, as above. Then $\lambda=0$, and we can assume that $f_3=h_3$. Combining (\ref{f3.16}) and Corollary  
 \ref{c1.2} we also conclude that   
\bes
a=\lambda w(f_1,f_2)+\mu, \ \ f=\a f_3+\lambda w(f_1,f_2)+\mu. 
\ees
If $\lambda\neq 0$, then $\tau$ admits a reduction of type IV. Otherwise $a=\mu\in F$ and we can apply Lemma \ref{l3.6}. 

\smallskip

Now, suppose that $\deg f_2\leq \deg f_3$. If $\overline{f_2}\in<\overline{f_1},\overline{f_3}>$, 
then $\overline{f_2}=\b \overline{f_3}+\gamma\overline{f_1^k}$. Hence Lemma \ref{l3.6} gives the triviality of (\ref{f3.7}). 

Let $\overline{b}=\overline{f_2}\notin<\overline{f_1},\overline{f_3}>$. 
Since $\deg\,b=\deg\,f_2\leq \deg\,f_3$, applying Corollary  
\ref{c1.1} gives $\deg\,f_3\geq N(f_1,f_3)$. By the definition of $N(f_1,f_3)$, we see that $f_1,f_3$ is a 2-reduced pair, i.e.  
$\deg\,f_1=2n$ and $\deg\,f_3=ns$, where $s\geq 3$ is an odd number.  
 According to Lemma \ref{l3.5}, we can assume that $g_2$ is an irreducible element of $\phi=(f_1,g_2,f_3)$. 
 Then $\phi$ satisfies the conditions of Proposition \ref{p1.2}. 

If $\deg_y(b(x,y))=k$, then (\ref{f1.1}) gives   
  $k=2$ if $s>3$, and $k=2,4$ if $s=3$. Consequently,  
$\deg_y(\frac{\d b}{\d y}(x,y))=k-1$ which equals $1$ or $3$, and (\ref{f1.1}) gives also 
$\deg(\frac{\d b}{\d y}(f_1,f_3))\geq ns$. We have  
\bee\label{f3.19}
[f_1,f_2]=[f_1,g_2]+[f_1,f_3]\frac{\d b}{\d y}(f_1,f_3).
\eee

Assume that Proposition \ref{p1.2}(1) is valid for $\phi$, i.e. 
$\deg g_2<n(s-2)+\deg[f_1,f_3]$. Comparing the degrees of elements in   (\ref{f3.19}) we obtain  
\bee\label{f3.20} 
\deg[f_1,f_2]\geq \deg[f_1,f_3]+ns. 
\eee
Consequently, $\deg f_2\geq n(s-2)+\deg[f_1,f_3]>n$ and $\overline{f_1}
\notin \langle\overline{f_2}\rangle$. Then either $\overline{f_1}$ and $\overline{f_2}$ are 
algebraically independent or $f_1, f_2$ is a $\ast$-reduced pair. 
Note that $N(f_1,f_2)>\deg\,[f_1,f_2]>ns\geq \deg\,a$. By Corollary \ref{c1.1}, in both cases we have  
$\overline{a}\in\langle\overline{f_1},\overline{f_2}\rangle$. 
Since  $\deg\,f_2^3>\deg\,f_3\geq \deg\,a$ and $\deg\,f_1f_2>\deg\,f_3$, it follows that   
\bes
a=\b f_2^2+\g f_2+G(f_1),\ \   \deg(G(f_1))\leq\deg f_3 ; 
\ees
moreover, we can have $\b\neq 0$ only if $\deg\,f_3\geq 2 \deg\,f_2$. Note that the equality \\
$\deg(G(f_1))=\deg\,f_3$ is impossible.  
Consequently, $\deg(G(f_1))<\deg\,f_3$. Since $f_1,f_3$ is a $2$-reduced pair, it follows that $f_1,\a f_3+G(f_1)$ is also a $2$-reduced pair. We have  
\bes 
\tau=(f_1,f_2,\a f_3+\b f_2^2+\g f_2+G(f_1)).
\ees
If $(\b,\g)\neq (0,0)$, then it is easy to check that $(f_1,g_2,\a f_3+G(f_1))$ is a reduction of 
 types I--III of the automorphism $\t$  with the active element $f_2$. Proposition \ref{p2.1} gives 
$d(\tau)=n(s+2)+\deg\,f_2=d=\deg\,\te$ and $t(\tau)\geq 1$. This contradicts the inequality $(d(\tau),t(\tau))\leq (d,t)$.  
If $\b=\g=0$, then we apply Lemma \ref{l3.6}.

\smallskip

We now consider the case when Proposition \ref{p1.2}(3) is valid for $\phi$. 
Then (\ref{f3.19}) again gives (\ref{f3.20}). Besides, in this case $s=3$ and $\deg f_2>\deg g_2=\frac{3n}{2}$. 
Consequently, the same argument as above gives $a=\b f_2+\g f_1+\delta$. By Lemma \ref{l3.6}, we can assume that $\b\neq 0$. 
Then  $(f_1,g_2,\a f_3+\g f_1+\delta)$ is a reduction of type 
I or II of $\t$ with the active element $f_2$. 
The value of $(d(\tau),t(\tau))$, which can be calculated by Proposition \ref{p2.1}, again contradicts the inequality $(d(\tau),t(\tau))\leq (d,t)$.

At last we consider the case when Proposition \ref{p2.1}(2) is valid for $\phi$.   
Let $(f_1,h_2,f_3-\delta h_2^2)$ be a reduction of type IV of $\phi$, where  
\bes
g_2=h_2+g,\ g\in\langle f_1,f_3\rangle\setminus F, \ \deg
(f_3-\delta h_2^2)\leq 2n.
\ees
Since $\deg g_2\leq \deg h_2  =\frac{3n}{2}$, we have $\deg g\leq
\frac{3n}{2}$. By Definition \ref{d1.4}, Proposition \ref{p1.2}(3) is valid for $(f_1,h_2,f_3)$.

If $\deg f_2>\frac{3n}{2}$, then by Lemma \ref{l3.5} we can assume that $\phi=(f_1,h_2,f_3)$ and this case reduces to the preceding one. 

Assume that $\deg f_2\leq \frac{3n}{2}$. By Proposition \ref{p2.1}, we have $d(\phi)=\frac{13n}{2}$, $t(\phi)=1$. 
Since $d=\deg\,\te\leq \frac{13n}{2}$, this contradicts the inequality $(d(\phi),t(\phi))<(d,t)$. 
 $\Box$

 \begin{lm}\label{l3.8}
If $f_1, f_3$ is a $\ast$-reduced pair and $\deg f_3 <\deg f_1$,
then (\ref{f3.7}) is trivial. 
\end{lm}
\Proof 
We first consider the case when $\deg f_1<\deg f_2$. Since 
$\deg\,a\leq \deg\,f_3<\deg f_1 <\deg f_2$, applying Lemma \ref{l3.6} we conclude that 
$\overline{a}\notin \langle\overline{f_1},\overline{f_2}\rangle$. Consequently, 
$\overline{f_1}$ and $\overline{f_2}$ are algebraically dependent. By Lemma \ref{l3.6},  
we assume that $\overline{f_2}\notin \langle\overline{f_1}\rangle$. Then $f_1,f_2$ is 
a $\ast$-reduced pair. Corollary \ref{c1.1} gives also $\deg\,a\geq N(f_1,f_2)$. 
Combining these inequalities and the definition of $N(f_1,f_2)$ we conclude that  
$\deg f_1=2n$ and $\deg f_2=3n$. Consider an elementary reduction $(f_1,f_2,g_3)$ of 
 $\te$ with an irreducible $g_3$ (assume that $g_3=f_3$  if $f_3$ is an 
irreducible element of $\te$). Then $(f_1,f_2,g_3)$ satisfies the conditions of 
Proposition \ref{p1.2}. By the same discussions related to the automorphism $(f_1,f_2,g_3)$ 
as in the proof of Lemma \ref{l3.7}, we obtain the triviality of  
 (\ref{f3.7}). 

\smallskip

If $\deg\,f_1=\deg\,f_2$, then by Lemma \ref{l3.6} we may assume that  $\overline{f_1}$ and $\overline{f_2}$ are linearly independent. 
Then $\overline{f_1}$ and $\overline{f_2}$ are algebraically independent. Consequently, $a\in F$. 

\smallskip

It remains to consider the case when $\deg\,f_2<\deg\,f_1$. Assume that  
$\overline{f_2}\in \langle\overline{f_3}\rangle$. By Lemma \ref{l3.6} we can assume that $\overline{f_2}$ and $\overline{f_3}$ are linearly independent.  
 Then $\deg\,f_2=t\cdot \deg\,f_3$, $t\geq 2$. 
 In this case $\overline{f_1}\notin \langle\overline{f_2}\rangle$, since $f_1,f_3$ is a $\ast$-reduced pair. 
 Consequently, $f_1,f_2$ is also a $\ast$-reduced pair. 
Since $\deg\,a\leq \deg\,f_3< \min\{\deg\,f_1,\deg\,f_2\}$, applying Corollary \ref{c1.1} and Lemma \ref{l3.6} we obtain  
  $\deg f_2=2n$, $\deg f_1=3n$, and 
\bes
\deg\,f_3\geq\deg\,a\geq N(f_1,f_2)=n+\deg\,[f_1,f_2]>\frac{2n}{t}=\deg f_3, 
\ees
which gives a contradiction.  

Therefore we can assume that $\overline{f_2}\notin \langle\overline{f_3}\rangle$, i.e.  
$\overline{b}=\overline{f_2}\notin \langle\overline{f_1},\overline{f_3}\rangle$. 
Then Corollary \ref{c1.1} gives $N(f_1,f_3)\leq \deg\,f_2<\deg\,f_1$. 
By the definition of $N(f_1,f_3)$ this is possible only if $f_1,f_3$ is a 2-reduced pair, 
and $\deg\,f_3=2n$, $\deg\,f_1=ns$, where $s\geq 3$ is an odd number. 
By Lemma \ref{l3.5} we can assume that $g_2$ is an irreducible element of $\phi=(f_1,g_2,f_3)$. 
Then $\phi$ satisfies the conditions of Proposition \ref{p1.2}. 
Applying the same part of discussions related  
to the automorphism $\phi$ as in the proof of Lemma \ref{l3.7}, we obtain the triviality of (\ref{f3.7}). $\Box$

\begin{lm}\label{l3.9}
If $\overline{f_1}\in \langle\overline{f_3}\rangle$ and $\deg f_1>\deg f_3$,
then the relation (\ref{f3.7}) is trivial. 
\end{lm}
\Proof 
Assume that $n=\deg\,f_3$ and $\overline{f_1}=\b \overline{f_3}^k$, $k\geq 2$. Then  $\deg f_1=nk$. By (\ref{f3.3}) we have  
\bes
\varphi_{i_1+t}=\s(2,1,b(x_1,x_2))=\s(1,1,-\b x_3^k)\s(2,1,b(x_1,x_2))\s(1,1,\b x_3^k).
\ees
After such replacement, instead of $\phi\rightarrow\te$ we obtain   
\bes
\phi\rightarrow (f_1-\b f_3^k,g_2,f_3)\rightarrow (f_1-\b f_3^k,f_2,f_3)\rightarrow \te. 
 \ees
Since $\deg\,(f_1-\b f_3^k,g_2,f_3),\deg\,(f_1-\b f_3^k,f_2,f_3)<d=\deg\,\te$, the new sequence (\ref{f3.8}) has the same exponent.  
Then instead of $\phi$ we can take $(f_1-\b f_3^k,f_2,f_3)$, i.e. we can assume that $\phi$ reduces 
the element $f_1$ of $\te$. Interchanging the elements $f_1$ and $f_2$ and applying Lemmas \ref{l3.6}, \ref{l3.7}, and 
\ref{l3.8} we can restrict ourselves to the case when $\overline{f_2}\in \langle\overline{f_3}\rangle$, $\deg\,f_2>\deg\,f_3$. 
Then $\deg f_2=nr$, $r\geq 2$. By Lemma \ref{l3.6} we may assume that $f_1,f_2$ is a $\ast$-reduced pair, 
i.e.   $\deg\,f_1=mp$, $\deg f_2=ms$, $(p,s)=1$, $m\geq n$. Hence   
$N(f_1,f_2)>m(ps-p-s)>n\geq \deg a$. Corollary \ref{c1.1} gives  
$\overline{a}\in\le\overline{f_1},\overline{f_2}\re$, i.e. $a\in F$. $\Box$

\smallskip

This finishes the proof of Theorem \ref{t3.1}.  

\bigskip

Now we formulate some problems related to the relations (\ref{f3.2})--(\ref{f3.4}) and discuss why it is important to study them. 

Let us denote by $G(A_n)$ the group defined by the system of generators (\ref{f3.1}) and the system of relations (\ref{f3.2})--(\ref{f3.4}). 

\begin{prob}\label{prob1}
Find a canonical form of elements of the group $G(A_n)$. In particular, is the word problem decidable in this group? 
\end{prob}

Note that the group $G(A_2)\cong GA_2(F)$ has 
a nice representation (see \cite{Cohn}). By Theorem \ref{t3.1}, if $F$ is a field of characteristic $0$, 
then $G(A_3)\cong TA_3(F)$. Note that the group $GA_n(F)=Aut\,A_n$ has the decidable word problem. In fact, if 
$\phi = (f_1,f_2,\ldots,f_n)\in GA_n(F)$, then $\phi = id$ iff $f_i=x_i,\, 1\leq i\leq n$.
Consequently, the groups   $G(A_2)$ and $G(A_3)$ have the decidable word problem. 
But we have no canonical form of elements of $G(A_3)$. 

Now, consider the homomorphism 
\bes
s_n : G(A_n) \longrightarrow TA_n(F) 
\ees
which sends the generators (\ref{f3.1}) to the corresponding automorphisms of $A_n$. 

\begin{prob}\label{prob2}
Find the kernel of $s_n$. 
\end{prob}

We know that the kernel of $s_2$ is trivial. By Theorem  \ref{t3.1} the kernel of $s_3$ is also trivial if $char(F)=0$.  
A solution of Problem \ref{prob2} gives a description of the group $TA_n(F)$ by generators and relations.

\bigskip

\begin{center}
{\bf\large Acknowledgments} 
\end{center}

\hspace*{\parindent}

I am grateful to I.\,Shestakov for thoroughly going over the details
of the proofs. I am also grateful to Max-Planck Institute f\"ur Mathematik for hospitality and exellent working conditions. I also thank J.\,Alev, P.\,Cohn, N.\,Dairbekov, M.\,Jibladze, L.\,Makar-Limanov, D.\,Wright, and M.\,Zaidenberg for numerous helpful comments and discussions. 

\bigskip

\hspace*{\parindent}

\end{document}